\documentclass[french,11pt,twoside]{amsart}

\usepackage{amsfonts}
\usepackage{amssymb}
\usepackage{amsmath}
\usepackage{latexsym}
\usepackage{amsbsy}
\usepackage{eufrak}
\usepackage[french]{babel}
\usepackage[latin1]{inputenc}
\usepackage[T1]{fontenc}
\usepackage{graphicx}
\usepackage{amscd}

\newtheorem{defn}{D\'efinition}[section]
\newtheorem{theoreme}[defn]{Th\'eor\`eme}
\newtheorem{proposition}[defn]{Proposition}
\newtheorem{lemma}[defn]{Lemme}
\newtheorem{corollary}[defn]{Corollaire}
\newtheorem{remark}[defn]{Remarque}

\numberwithin{equation}{section}

\theoremstyle{plain}

\newtheorem{example}{Exemple}

\numberwithin{equation}{section}

\def\preuve{\smallskip\goodbreak{\it Preuve.~~--~\kern.3em}
     \ignorespaces}%
\def\qedbox{$\square$}%
\def\qed{\ifmmode\qedbox\else\unskip\ \hglue0mm\hfill
     \qedbox\smallskip\goodbreak\fi}%
\def\preuvedt#1{\smallskip\goodbreak{\it Preuve du th\'eor\`eme~\ref{#1}.~~--~\kern.3em}
     \ignorespaces}%
\def\preuvedp#1{\smallskip\goodbreak{\it Preuve de la proposition~\ref{#1}.~~--~\kern.3em}
     \ignorespaces}%
\def\preuvedl#1{\smallskip\goodbreak{\it Preuve du lemme~\ref{#1}.~~--~\kern.3em}
     \ignorespaces}%
\newcommand{\eps}{\varepsilon}
\newcommand{\equiva}{\Longleftrightarrow}
\newcommand{\Om}{\Omega}
\newcommand{\Vol}{\mathrm{Vol}}
\newcommand{\ii}{\iota}

\newcommand{\N}{\mathbb N}

\newcommand{\R}{\mathbb R}
\newcommand{\C}{\mathbb C}

\begin{document}

\title[]{A mi-chemin entre analyse complexe et superanalyse.}

\maketitle
\author{Pierre Bonneau, Anne  Cumenge}\\
\author{}
\address{Equipe Emile Picard, Institut de Mathématiques, {\sc umr} 5219\\
 Universit\'{e} Paul Sabatier - 31062 Toulouse Cedex 9}\\

\email{{\it pierre.bonneau@math.univ-toulouse.fr}, {\it
anne.cumenge@math.univ-toulouse.fr}}

 \setcounter{part}{-1}

\smallbreak \noindent {\bf Abstract.}

\noindent In the framework of superanalysis we get a functions theory close to complex analysis,
under a suitable condition (A) on the real superalgebras in consideration (this condition is a generalization
of the classical relation $1 + i^2 = 0$ in $\C$). Under the condition (A), we get an integral representation formula
for the superdifferentiable functions. We deduce properties  of  the superdifferentiable functions : analyticity,
a result of separated superdifferentiability, a Liouville theorem and a continuation theorem of Hartogs-Bochner type.

\smallbreak \noindent {\bf Keywords : } superspace ; superdifferentiable function ; integral representation ; superanalysis.

\smallbreak \noindent 30E20 - 32A26- 30G30- 35C15

\setcounter{section}{-1}
\section{Introduction.}\setcounter{section}{0}

\smallbreak L'analyse complexe s'est fortement développée au 19°
siècle, en particulier avec les travaux de Cauchy, Riemann et
Weierstrass. Les essais de généralisation en dimension supérieure
conduisirent à la construction des quaternions et des octonions.
Toutefois, ces nouveaux objets s'avérèrent alors un peu décevants
au niveau de l'analyse. Par ailleurs, de fortes limitations à ce
type de généralisations furent obtenues, entre autres par
Frobenius avec une classification des algèbres associatives
réelles de division de dimensions finies, ou par Bott-Milnor et Kervaire avec le résultat suivant :\\
\noindent{Théorème (\cite{BM}, \cite{Ke}):}  L'espace vectoriel
$\mathbb{R}^{n}$ possède une opération produit $\R$-bilinéaire
sans diviseur de 0 seulement pour
n=1,2,4 ou 8. \\
 Pour n=4 et 8, on retrouve, entre autres,   les
quaternions et les octonions. \\
Nous considèrerons seulement  dans toute la suite des algèbres
associatives. Les algèbres associatives de division étant
unitaires,    la condition d'intégrité  dans le cas de
$\mathbb{R}^{2}$ est équivalente à l'existence d'un neutre
multiplicatif 1 et d'un élément i de $\mathbb{R}^{2}$ tel que
$1+i^{2}=0$.

Sachant qu'il est impossible de trouver en dimensions strictement
supérieures à 4 des $\R$-algèbres associatives de division, nous
rechercherons des algèbres réelles (associatives) ayant de bonnes
propriétés analytiques, c'est-à-dire dans lesquelles existe une
théorie des fonctions analogue à l'analyse complexe. Nous
essaierons d'étendre, en dimension supérieure à deux, le point de
vue de Cauchy-Riemann sur l'analyse
complexe (fonctions dérivables et représentation intégrale).\\

Des conditions de Cauchy-Riemann apparaissent naturellement en
superanalyse (cf. \cite{K}, \cite{R} par exemple). L'analyse  sur les superespaces,
dite "superanalyse", s'est développée dans les décennies '60 et '70, formalisant les superespaces des
physiciens, espaces dans lesquels cohabitent des variables anti-commutatives avec les variables usuelles.
Pour un historique et des références, nous renvoyons aux livres de F.A. Berezin \cite{Be}, A. Khrennikov \cite{K} et A. Rogers \cite{R}.\\
Nos résultats sont énoncés dans le cadre des  superespaces.
Soulignons qu'en provenance   de la superanalyse, nous utiliserons
seulement les définitions de superespaces et de la
superdifférentiabilité ; nos méthodes sont inspirées par celles de
l'analyse complexe.

\textit{}Il est à noter que, contrairement à l'analyse complexe,
quaternionique ou de Clifford, nous ne disposons pas d'un
opérateur de conjugaison sur les algèbres $\Lambda=\Lambda_0
\oplus \Lambda_1$ -- où les éléments de $\Lambda_0$ (resp.
$\Lambda_1$) commutent (resp. anti-commutent) entre eux -- et les
superespaces $\Lambda_0^{\, n}\times \Lambda_1^{\, m}$ associés à
$\Lambda$ que nous considérons.

\smallbreak Nous définirons sur certains  superespaces  un
opérateur $d''$ de Cauchy-Riemann dont le noyau coïncidera avec
l'espace des fonctions superdifférentiables  (S-différentiables en abrégé). \\
Des conditions ($A_j$), $j=0,1$ sur $\Lambda_j$, à mettre en
parallèle avec la relation complexe $1+i^2 =0$, s'imposent alors
de manière naturelle lorsque nous cherchons une solution
fondamentale de cet opérateur $d''$. Plus précisément : \\
\noindent $\quad\rm{(A_{0})}$ \;   il existe une base
$(e_{0}=1,e_{1},...,e_{p})$  de $\Lambda_{0}$ vérifiant
$\displaystyle{\sum_{k=0}^{p}\,e_{k}^{\;2}=0}$ \,,

\noindent $\rm{(A_{1})}$ $\begin{array}{l} {\mbox{ il existe une
base }} (\eps_{1},...,\eps_{q} ) {\mbox{ de }} \Lambda_{1} {\mbox{
et une suite finie }}
s_{1}=1 <s_{2}<...<s_{r}<s_{r+1}=q+1 \\
{\mbox{  telles que, pour tout }} j=1,...,q,
 {\mbox{ il existe }} a_{j}\in\Lambda_{0} {\mbox{ vérifiant }}
\eps_{j}=a_{j}\eps_{s_{k}}\, {\mbox{ si }} \,s_{k}\leq j<s_{k+1},\\
{\mbox{  avec }}
a_{s_{1}}=a_{s_{2}}=...=a_{s_{r}}=e_{0}\; {\mbox{ et }}\;
\sum_{j=s_{k}}^{s_{k+1}-1}\,a_{j}^{\;2}=0\; {\mbox{ pour tout }}
k=1,...,r\,.
\end{array}$

\smallbreak\noindent La condition $(A_0 )$ est une condition
algébrique nécessaire (et suffisante) pour l'obtention d'une
solution fondamentale de notre opérateur $d''$ opérant sur les
formes définies sur $\Lambda_0$ ; quant à la condition $(A_1 )$,
elle est suffisante et "presque" nécessaire si l'on espère un
caractère explicite pour une solution fondamentale de l'opérateur
$d''$ sur $\Lambda_1$.

 \smallbreak  Nous obtenons, à partir d'une solution fondamentale pour l'opérateur $d''$,
une formule de représentation des formes différentielles avec
opérateurs intégraux à noyaux $K(x,y)$ explicites ; et en
particulier, pour les fonctions :
\begin{theoreme} On suppose les conditions ($A_0$) et ($A_1$) satisfaites par l'algèbre $\Lambda= \Lambda_0 \oplus
\Lambda_1$. Soit $D$ ouvert de
$\Lambda_{0}^{n}\times\Lambda_{1}^{m}$ borné et à frontière $C^1$
et $f$ fonction  de classe $C^{1}$ au sens de Fréchet dans $D$,
continue sur $\overline{D}$ ainsi que $df$, alors pour tout $x\in
D\;$ :
\begin{eqnarray*}
 f(x)= \int_{\partial D}\, f(y) K(y,x) \,- \,\int_{D}\, d''f(y) \, {\scriptstyle{\wedge}} K(y,x)
\; .
\end{eqnarray*}
\end{theoreme}

 \smallbreak
\noindent Nous déduisons de la formule de représentation des
fonctions qS-différentiables (i.e. vérifiant $d''f=0$) des
propriétés de ces fonctions, comme l'harmonicité, une propriété de
qS-analyticité, un théorème de Liouville, ainsi qu'un résultat de
qS-différentiabilité séparée  :

\begin{theoreme} :  Si $f$ est une
fonction  définie sur un domaine $D$ de $\Lambda_0^{\,n} \times
\Lambda_1 ^{\,m}$ à valeurs dans $\Lambda$ séparément
qS-différentiable sur $D$ par rapport à chacune de ses
"hypervariables" appartenant à $\Lambda_0$ ou $\Lambda_1$, alors
$f$
 est qS-différentiable sur $D$.
\end{theoreme}

et un théorème de prolongement  de type Hartogs :

\begin{theoreme}: Sous les conditions ($A_0$) et ($A_1$), si $\partial \Om$ est le bord connexe d'un domaine
$\Om$ borné de $\Lambda_0 ^{\,n} \times \Lambda_1 ^{\, m}$, avec
$n+m \geq 2$,  et $f$ une fonction qS-différentiable dans un
voisinage connexe de $\partial \Om$, alors $f$ se prolonge en une
fonction qS-différentiable sur $\Om$.
\end{theoreme}

\medbreak Dans une première partie, nous préciserons nos notations
et rappellerons les principales notions de superanalyse.
Nous fournirons aussi quelques exemples de superalgèbres. \\
Dans la seconde partie, suivant le point de vue "riemannien" de
l'analyse complexe, nous définirons un opérateur de Cauchy-Riemann
dont le noyau est constitué des fonctions
super-différentiables en les variables commutatives et "quasi"-super différentiables en les variables anti-commutatives. \\
La troisième partie sera dévolue à la recherche de représentations
intégrales pour les fonctions définies dans un super-espace et à
valeurs dans notre algèbre. C'est dans ce paragraphe
qu'apparaîtront les conditions $(A_0$) et ($A_1$), qui fourniront
une superanalyse aux propriétés étonnamment voisines de l'analyse
complexe. Les représentations intégrales seront alors l'outil que
nous utiliserons pour étudier les propriétés des fonctions
super-différentiables, ce qui fait l'objet de la quatrième partie.
Dans une dernière partie, nous regroupons quelques commentaires
sur les conditions algébriques $(A_0 )$ et $(A_1 )$.\\
L'objectif
de l'article n'est pas un recensement exhaustif des propriétés des
fonctions super-différentiables ; il essaie de souligner une
étrange proximité, sous certaines conditions, entre l'analyse
complexe et la superanalyse. Ces résultats ont été partiellement annoncés dans \cite{BC}.

\medbreak  Les auteurs tiennent à remercier le referee pour ses intéressants commentaires.

\section{Superespaces.}  Nous noterons $\Lambda$ une
{\it superalgèbre commutative réelle (en abrégé CSA)}, dont nous rappelons ci-dessous la définition (cf.  \cite{K} ou \cite{R} par exemple). \\
Un
$\mathbb{R}$-espace vectoriel $Z_2$-gradué $\Lambda=\Lambda_{0}\bigoplus\Lambda_{1}$  muni d'une fonction parité $\sigma$ définie sur les éléments homogènes par $\sigma (a)=0$ (resp. 1) si $a\in \Lambda_0$ (resp. si $a\in \Lambda_1$) devient une superalgèbre si on le munit d'une structure d'algèbre associative unitaire dont la multiplication vérifie la propriété  $\sigma (ab)= \sigma (a) +\sigma (b)\;\mbox{(mod. 2)} $ pour tout couple $(a,b)$ d'éléments homogènes.\\
Le supercommutateur est défini sur les éléments homogènes par $[a,b\}=ab -(-1)^{\sigma (a)\sigma(b)} ba$ (et prolongé par bilinéarité).
Une superalgèbre est dite commutative si pour tout couple $(a,b)$ d'éléments homogènes $[a,b\}=0$.\\
Nous définissons le $\Lambda_{1}$-annihilateur comme étant
$^{\bot}\Lambda_{1}=\{\lambda\in\Lambda :\lambda\Lambda_{1}=0\}.$
Toutes les CSA considérées ici seront de dimensions finies. Nous
noterons $(e_{0}, e_{1}, ...,e_{p})$ une base de $\Lambda_{0}$
($e_{0}$ est l'élément unité de $\Lambda$), $(\eps_{1},\eps_{2},
...,\eps_{q})=(e_{p+1},e_{p+2}, ...,e_{p+q})$ une base de
$\Lambda_{1}$, et définissons les coefficients de structure
(d'algèbre) $\Gamma$ de $\Lambda$ par
$e_{i}e_{j}=\sum_{k=0}^{p+q}\,\Gamma_{i,j}^{k}e_{k}$ pour
$i,j=0,1, ...,p+q$. On remarque que, d'après la définition d'une
CSA, les coefficients $\Gamma$ sont symétriques en $(i,j)$ si
$i\,ou\,j\,\in\{0,1,...,p\}$, anti-symétriques si
$i\,et\,j\in\{p+1,p+2,...,p+q\}$. \\

On appelle {\it superespace sur la CSA $\Lambda$} tout
$\mathbb{R}$-espace vectoriel
$$\mathbb{R}_{\Lambda}^{n,m}\,=\Lambda_{0}^{n}\times\Lambda_{1}^{m}$$
 où $n,m\in\mathbb{N}$. \\
\begin{defn}\label{S-diff}
Si $U$ est un ouvert de $\mathbb{R}_{\Lambda}^{n,m}$ et $F$ une
application de $U$ dans $\Lambda$, nous disons que {\it $F$ est
superdifférentiable à droite} (ou {\it S-différentiable}) en $x\in  U$ s'il
existe des éléments $\frac{\partial F}{\partial x_{j}}(x)$ de $\Lambda$,
$j=1,...,n+m$ tels que, pour tout $h\in\mathbb{R}_{\Lambda}^{n,m}$
tel que $x+h\in U$, on ait
$$F(x+h)=F(x)\,+\,\sum_{j=1}^{n+m}\frac{\partial F}{\partial x_{j}}(x)h_{j}\,+\,o(h)$$
\end{defn}

\noindent
avec $lim_{\|h\|\mapsto 0}\frac{\|o(h)\|}{\|h\|}=0$ où $\|\,.\,\|$
est une norme sur l'espace vectoriel $\mathbb{R}_{\Lambda}^{n,m}.$
\\
Nous remarquons que $\frac{\partial F}{\partial
x_{1}}(x),...,\frac{\partial F}{\partial x_{n}}(x)$ sont définis
de façon unique par la condition ci-dessus, tandis que
$\frac{\partial F}{\partial x_{n+1}}(x),...,\frac{\partial
F}{\partial x_{n+m}}(x)$ sont définis modulo $^{\bot}\Lambda_{1}.$
La condition de S-différentiabilité de $F$ exprime donc le fait
que $F$ est différentiable et que sa dérivée est définie par les
opérateurs de multiplication par des éléments de $\Lambda$.

\smallbreak \noindent \begin{example}\label{exemple-C} \textup{ Si
$\Lambda_{1}=\{0\}$  et $\Lambda_{0}= \mbox{ Vect }(e_0 ,e_1 )$
avec $e_{0}=1$ et $e_{1} ^{\, 2}=- e_0$, les fonctions
S-différentiables sont les fonctions holomorphes et la
superanalyse est alors l'{\it analyse complexe}.}
\end{example}

\smallbreak \noindent

\begin{example}\label{exemple-hyperbolique} ({\it Analyse
hyperbolique.}) $\Lambda_{1}=\{0\}$ et $\Lambda_{0}= \mbox{ Vect
}(e_0 , e_1)$ \textup{ avec
$e_{0}=1$ et $(e_1 )^{2}=e_0$. \\
$f= u\,e_0 +v\, e_1$ $\mathcal{S}$-diff. $\equiva$ $\frac{\partial
 u}{\partial x_0 }=\frac{\partial
 v}{\partial x_1 }$ et $\frac{\partial
 u}{\partial x_1 }=\frac{\partial
 v}{\partial x_0 }$  \\
 Alors $f$, $u$ et $v$ vérifient l'équation des ondes }
 $\frac{\partial ^2
 u}{\partial x_0 ^{\,2}} - \frac{\partial ^2
 u}{\partial x_1 ^{\,2 }}=0$.
\end{example}

 \smallbreak\noindent Remarque : soient $\varphi$ et $\psi$
fonctions de classe $C^{1}$ sur $\mathbb{R}.$ Alors
$f:\;\Lambda\rightarrow\Lambda$ définie par \\
$f(xe_{0}+ye_{1})=[\varphi(x+y)+\psi(x-y)]e_{0}
+[\varphi(x+y)-\psi(x-y)]e_{1}$ est S-différentiable, mais pas
plus régulière que ne le sont $\varphi$ et $\psi$.

\smallbreak\noindent

\smallbreak
\begin{example}\label{exemple-table}
 \textup{Supposons $\Lambda_{1}$ et
$\Lambda_{0}$ de dimensions  6, avec la table de multiplication
suivante :} \smallbreak
\begin{center}
\begin{tabular}{|c||c|c|c|c|c|c|c|c|c|c|c|c|}

 \hline

 $\times$ & $\, e_0 \,$ & $e_1$   & $\, e_2 \,$ & $e_3$  & $e_4$ &$e_5 $ & $\eps_1$   & $\eps_2$& $\eps_3$& $\eps_4$&$\eps_5$&$\eps_6$ \\
 \hline
 \hline
 $e_0$    & $e_0$ & $e_1$   & $e_2$ & $e_3$ & $e_4$ & $e_5$  & $\eps_1$   & $\eps_2$& $\eps_3$& $\eps_4$ &$\eps_5$&$\eps_6$\\
 \hline
 $e_1$ &    $e_1$ & $- e_0$   & $e_3$ & $- e_2$ &$e_5$ & $- e_4$  & $\eps_2$   & $- \eps_1$& $\eps_6$& $ \eps_5$&$-\eps_4$&$-\eps_3$ \\
 \hline
 $e_2$ & $e_2$ & $e_3 $ &0&0&0&0&0 &0 &0&0&0&0\\
 \hline
 $e_3$ & $e_3$ & $-e_2 $ &0&0&0&0&0&0 &0 &0&0&0\\
 \hline
$e_4$ & $e_4$ & $e_5 $ &0&0&$e_2$ &$e_3$&$\eps_3$ &$\eps_2 $&0&$\eps_3$& $\eps_6$&0\\
 \hline
 $e_5$ & $e_5$ & $- e_4 $ &0&0&$e_3$ &$-e_2$&$\eps_2$&$-\eps_3$ &0 &$\eps_6$&$-\eps_3$&0\\
 \hline
 $\eps_1 $ & $\eps_1 $ & $\eps_2 $ &0&0&$\eps_3$&$\eps_2$&0&0&0& $e_2 $ & $e_3 $ &0\\
 \hline
 $\eps_2 $ & $\eps_2 $ & $-\eps_1 $ &0&0&$\eps_2$&$-\eps_3$&0&0&0& $e_3 $ & $-e_2 $&0 \\
 \hline
 $\eps_3 $ & $\eps_3 $ & $\eps_6 $ &0&0 &0&0& 0 & 0 &0&0 &0&0 \\
 \hline
 $\eps_4 $ & $\eps_4 $ & $\eps_5 $&0&0 &$\eps_3$&$\eps_6 $ & $-e_2 $ & $-e_3 $ &0&0&0&0 \\
 \hline
 $\eps_5 $ & $\eps_5 $ & $-\eps_4 $ &0&0&$\eps_6$&$-\eps_3$&$-e_3 $ & $e_2 $&0 &0&0&0\\
 \hline
 $\eps_6 $ & $\eps_6 $ & $-\eps_3 $ &0&0&0&0&0&0&0&0&0&0 \\
 \hline
 \end{tabular}
\end{center}
\end{example}

\medbreak

\noindent L'ensemble des nilpotents de $\Lambda_0$ étant ici Vect
($e_2 ,e_3 , e_4 , e_5 $), l'algèbre commutative $\Lambda_0$ n'est
pas semi-simple  et donc n'est pas isomorphe à un produit
d'algèbres $\prod_{j=1}^s L_j $, où  $L_j =\R$ ou $\C$.

\section{L'opérateur de
Cauchy-Riemann en superanalyse.}

Soit $f$ une fonction d'un ouvert $U$ de
$\mathbb{R}_{\Lambda}^{n,m}=\Lambda_{0}^{n}\times\Lambda_{1}^{m}$
dans $\Lambda=\Lambda_{0} \oplus \Lambda_{1}$ ;
$f(x)=\sum_{k=0}^{p}\,f_{k}(x)e_{k}\,+\,\sum_{l=1}^{q}\,f_{p+l}(x)\eps_{l}$,
où $f_{0},..., f_{p+q}$ sont des fonctions réelles.\\
 Il sera
commode, comme dans \cite{K}, d'écrire
$x\in\mathbb{R}_{\Lambda}^{n,m}$ sous la forme $x=(y,\theta)$ avec
$y=(y_1,\dots,y_n )\in\Lambda_{0}^{n}$ et $\theta=(\theta_1 ,\dots
,\theta_m )\in\Lambda_{1}^{m}$. Ainsi avec
$y_{i}=\sum_{k=0}^{p}\,y_{i}^{k}e_{k}\in\Lambda_{0}$ et
$\theta_{j}=\sum_{l=1}^{q}\,\theta_{j}^{l}\eps_{l}$ :

\begin{eqnarray*}
df(x)=\sum_{i=1}^{n}\sum_{k=0}^{p}\frac{\partial f}{\partial y_{i}^{k}}(x)dy_{i}^{k}\;+\;\sum_{j=1}^{m}
\sum_{l=1}^{q}\frac{\partial f}{\partial \theta_{j}^{l}}(x)d\theta_{j}^{l}
\end{eqnarray*}
Soit  $f$ S-différentiable en $x$ (voir la définition
\ref{S-diff}) et  $h=(h_{1},..., h_{n}, h_{1}',...,
h_{m}')\in\Lambda_{0}^{n}\times\Lambda_{1}^{m}$ ; nous avons
\begin{eqnarray*}
df(x)(h)&=&\sum_{i=1}^{n}\,\frac{\partial f}{\partial y_{i}}(x)h_{i}\,+\,\sum_{j=1}^{m}\,\frac{\partial f}{\partial \theta_{j}}(x)h_{j}' \\
        &=&\sum_{i=1}^{n}\sum_{k=0}^{p}\,\frac{\partial f}{\partial y_{i}}(x)\, e_{k}dy_{i}^{k}(h)\,
        +\,\sum_{j=1}^{m}\sum_{l=1}^{q}\,\frac{\partial f}{\partial \theta_{j}}(x)\,\eps_{l}d\theta_{j}^{l}(h)
\end{eqnarray*}
par conséquent
\begin{eqnarray*}
\forall i=1,...,n\;,\; \forall k=0,...,p \;,\;\frac{\partial f}{\partial y_{i}}(x)e_{k}=\frac{\partial f}{\partial y_{i}^{k}}(x),\\
\forall j=1,...,m \;,\; \forall l=1,...,q\;, \;\frac{\partial  f}{\partial\theta_{j}}(x)\eps_{l}=\frac{\partial  f}{\partial\theta_{j}^{l}}(x),
\end{eqnarray*}
ou encore, puisque $ \frac{\partial f}{\partial
 y_{i}}(x)=\frac{\partial  f}{\partial  y_{i}^{0}}(x)$ :\\

 \noindent {\sc propriété} : {\it une fonction $f$   Fréchet-différentiable en ses variables réelles est S-différentiable
 en $x \in \Lambda$ si et seulement si}
\begin{eqnarray}
&& \frac{\partial  f}{\partial  y_{i}^{k}}(x)=\frac{\partial  f}{\partial y_{i}^{0}}(x)e_{k}\;,\;\forall  i=1,...,n
\mbox{ et }\forall k=1,...,p \\
&& (\frac{\partial  f}{\partial\theta_{j}^{1}}(x),..., \frac{\partial  f}{\partial\theta_{j}^{q}}(x))\in
\Lambda(\eps_{1},...,\eps_{q}):= \{(a\eps_{1},...,a\eps_{q})\;,\, a \in \Lambda \}\;, \forall  j=1,...,m
\end{eqnarray}
\noindent Rappelons que $\frac{\partial
f}{\partial\theta_{j}^{l}}(x)$ est défini modulo
$^{\bot}\Lambda_{1}.$

Nous sommes ainsi amenés naturellement à traiter différemment les
"variables commutatives" $y_{1},\dots, y_{n}$ et les "variables
anti-commutatives" $\theta_{1},\dots , \theta_{m}$.

\subsection{Etude dans $\Lambda_0 ^n$.} Oublions provisoirement
les variables anti-commutatives et intéressons-nous aux variables
$y_{1},..., y_{n}$. Et, d'abord, considérons le cas n=1.

Soit $f$ une fonction d'un ouvert $U$ de $\Lambda_{0}$ dans
$\Lambda$. Si $y=\sum_{k=0}^{p}\,y^{k}e_{k}\in U$, alors $f$ est
S-différentiable en y si $df(y)=\frac{\partial f}{\partial
y^{0}}(y) \,\sum_{k=0}^{p}e_{k}dy^{k}$ puisque, pour tout k,
$\frac{\partial f}{\partial y^{k}}(y)=\frac{\partial f}{\partial
y^{0}}(y)e_{k}$ comme nous venons de le voir. Dans le cas général,
nous décomposerons $df(y)$ en une partie $\Lambda_{0}$-linéaire
notée $d'f(y)$ et un reste noté $d''f(y)$: $$df(y)=d'f(y) +
d''f(y)$$ avec $$d'f(y)=\frac{\partial f}{\partial
y^{0}}(y)\sum_{k=0}^{p}\,e_{k}dy^{k},$$
$$d''f(y)=\sum_{k=1}^{p}(\frac{\partial f}{\partial
y^{k}}(y)-e_{k}\frac{\partial f}{\partial y^{o}}(y))dy^{k}.$$
Plus
généralement, si $\omega$ est une forme différentielle définie
dans $U$, à valeurs dans $\Lambda$, nous définissons
\begin{eqnarray*}
d'\omega=\sum_{k=0}^{p}\,e_{k}dy^{k}\frac{\partial \omega}{\partial y^{0}}:=dY \frac{\partial \omega}{\partial y^{0}},\\
d''\omega=\sum_{k=1}^{p}\,dy^{k}(\frac{\partial \omega}{\partial y^{k}}-e_{k}\frac{\partial \omega}{\partial y^{0}}).
\end{eqnarray*}
On vérifie immédiatement que $d'^{2}=d''^{2}=0$ et donc
$d''d'+d'd''=0$. Une fonction f est S-différentiable dans U si et
seulement si $d''f=0$ dans U. On généralise aisément au cas de
plusieurs "variables commutatives". Si
$y=(y_{1},y_{2}...,y_{n})\in\Lambda_{0}^{n}$ avec
$y_{i}=\sum_{k=0}^{p}\,y_{i}^{k}e_{k}$, on définit pour une forme
différentielle $\omega$ de degré quelconque définie dans un ouvert
U de $\Lambda_{0}^{n}$, et à valeurs dans $\Lambda,$
\begin{eqnarray}\label{d''lambda0^n}
d''\omega=d''_y \omega=\sum_{i=1}^{n}\sum_{k=1}^{p}\,dy_{i}^{k}(\frac{\partial \omega}{\partial y_{i}^{k}}-e_{k}
\frac{\partial \omega}{\partial y_{i}^{0}}),\\
d'\omega=\sum_{i=1}^{n}\sum_{k=0}^{p}\,e_{k}dy_{i}^{k}\frac{\partial \omega}{\partial y_{i}^{0}} :=
\sum_{i=1}^{n}dY_i \frac{\partial \omega}{\partial y_{i}^{0}}\,.\nonumber
\end{eqnarray}

\noindent{\sc Remarque}: dans le cas de $\Lambda_{0}=\mathbb{C}$,
la manière traditionnelle de définir l'opérateur de Cauchy-Riemann
n'est pas la même que ci-dessus, puisque l'habitude est de faire
éclater df en une partie $\mathbb{C}$-linéaire et une partie
$\mathbb{C}$-antilinéaire. Cette procédure ne peut opérer ici car
nous n'avons pas, dans le cas général, d'opérateur de conjugaison.
Toutefois, l'opérateur de Cauchy-Riemann habituel dans
$\mathbb{C}$ et celui qui est défini ci-dessus ont le même noyau,
l'ensemble des fonctions holomorphes. Et c'est ce qui nous
importe, puisque notre objectif est l'étude des fonctions
S-différentiables.

\subsection{Etude dans $\Lambda_1 ^m$}Occupons nous, maintenant, des "variables non-commutatives", et,
tout d'abord, du cas d'une seule variable anti-commutative. Soit
$f$ une fonction définie dans un ouvert $U$ de $\Lambda_{1}$, à
valeurs dans $\Lambda$. Si
$\theta=\sum_{l=1}^{q}\,\theta^{l}\eps_{l}\in  U$, $f$ est
S-différentiable en $\theta$ s'il existe un élément
$\frac{\partial  f}{\partial\theta}(\theta)$ de $\Lambda$ tel que,
pour tout $l=1,...,q$, $\frac{\partial
f}{\partial\theta^{l}}(\theta) = \frac{\partial f}{\partial\theta
}(\theta)\eps_{l}.$

Mais, contrairement au cas des "variables commutatives", il
n'existe, en général, aucun moyen canonique pour déterminer le
coefficient multiplicatif $\frac{\partial
f}{\partial\theta}(\theta)$ de l'application
$\Lambda_{1}$-linéaire $d'f(\theta)$; autrement dit, il n'existe,
en général, aucun moyen naturel pour décomposer $df(\theta)$ en
une partie $\Lambda_{1}$-linéaire et un reste. Pour conserver le
maximum de souplesse dans le choix de ces projections associées de
$df(\theta)$, nous supposerons que nous disposons d'une norme
$\|\;\|$ dans $\Lambda$ (norme telle que $\|xx'\|\leq
C\|x\|\|x'\|$) déduite d'un produit scalaire g. \\

Nous cherchons à définir un opérateur $d''$ dont le noyau contient
les fonctions S-différentiables.\\
 Nous pourrions effectuer une
projection orthogonale de $\left(\frac{\partial f}{\partial \theta ^1 }
(\theta), \dots ,\frac{\partial f}{\partial \theta ^q}(\theta)\right)$
sur $\Lambda (\eps_1 ,\dots \eps_q )$ afin de définir d'abord $d'
=d - d''$ puis $d''$, mais cela ne donne pas, en général, un
opérateur $d''$ vraiment explicite et, partant, il semble
difficile d'obtenir ainsi une solution fondamentale explicite pour
l'opérateur $d''$, ce qui
sera l'objectif du paragraphe suivant. \\
$d''$ est bien plus explicite si nous supposons satisfaite la
condition suivante :
 \begin{eqnarray}\label{préA1}
  & &\mbox{ \it il existe une suite finie } 1=s_{1}<s_{2}<...<s_{r}<s_{r+1}=q+1
 \;\mbox{ \it telle que } \nonumber \\
 && \eps_{j}=\eps_{s_{k}}a_{j}\;si\;s_{k}\leq j<s_{k+1} \,\mbox{ \it où } a_{j}\in \Lambda
 \; \mbox{ \it et }\;
  a_{s_{1}}=a_{s_{2}}=...=a_{s_{r}}=e_{0}
\end{eqnarray}
  (en fait $a_{j}$ est un élément de
 $\Lambda_{0}$).\\
 Alors :\\

 \noindent {\sc propriété} : {\it Sous la condition (\ref{préA1}) sur $\Lambda_1$,
 une fonction $f$ Fréchet-différentiable en $\theta$ est S-différentiable en
 $\theta \in \Lambda_1$ si}
\begin{equation*}
 (\frac{\partial f}{\partial\theta^{s_{1}}}(\theta ), \dots ,\frac{\partial
f}{\partial\theta^{s_{r}}} (\theta))\in \Lambda (\eps_{s_{1}},\dots ,\eps_{s_{r}} )
\end{equation*}
\begin{equation}\label{qSdiff-1}
\frac{\partial  f}{\partial\theta^{j}}(\theta)=\frac{\partial  f}{\partial\theta^{s_{k}}}(\theta)a_{j} \;
\mbox{ si }\;s_{k}\leq j<s_{k+1}
\end{equation}
La plus utile des deux relations ci-dessus lorsque nous
chercherons à donner des propriétés de régularité des fonctions
S-différentiables sur un ouvert de $\Lambda_1$ sera
(\ref{qSdiff-1}). Nous sommes donc amenés à poser
$$d''(=d_{\theta}'')=\sum_{k=1}^{r}\sum_{j=s_{k}+1}^{s_{k+1}-1}d\theta^{j}\left(\frac{\partial}{\partial\theta^{j}}-a_{j}\frac{\partial}{\partial\theta^{s_{k}}}\right),$$
et
$$d'=d- d'' =\sum_{k=1}^{r}\sum_{j=s_{k}}^{s_{k+1}-1}d\theta^{j}a_{j}\frac{\partial}{\partial\theta^{s_{k}}}$$

On vérifie immédiatement que $d'^{2}=d''^{2}=0$ et $d'd''+d''d'=0.$\\
On généralise au cas de plusieurs "variables anti-commutatives":
Si
$\theta=(\theta_{1},\theta_{2},...,\theta_{m})\in\Lambda_{1}^{m}$
avec $\theta_{i}=\sum_{j=1}^{q}\,\theta_{i}^{j}\eps_{j}$, on
définit pour une forme différentielle $\omega$ définie dans un
ouvert $U$ de $\Lambda_{1}^{m}$,

\begin{equation}\label{d"lambda1^m}
\begin{array}{l}
d'\omega(=d_{\theta}'\omega)=\displaystyle{\sum_{i=1}^{m}\sum_{k=1}^{r}\sum_{j=s_{k}}^{s_{k+1}-1}\,a_{j}d\theta_{i}^{j}
 \frac{\partial\omega}{\partial\theta_{i}^{s_{k}}}
 {\mbox{ \quad noté }} \,d'\omega=\sum_{i=1}^{m}\sum_{k=1}^{r} dZ_k (\theta_i )
 \frac{\partial\omega}{\partial\theta_{i}^{s_{k}}}} \\
d''\omega(=d_{\theta}'' \omega)=\displaystyle{\sum_{i=1}^{m}\sum_{k=1}^{r}\sum_{j=s_{k}+1}^{s_{k+1}-1}\,d\theta_{i}^{j}
(\frac{\partial\omega}{\partial\theta_{i}^{j}}-\frac{\partial\omega}{\partial\theta_{i}^{s_{k}}}a_{j})}
\end{array}
\end{equation}

\subsection{Cas général} Enfin, dans le cas d'un superespace $\mathbb{R}_{\Lambda}^{n,m}
 =\Lambda_{0}^{n}\times\Lambda_{1}^{m}$  sur la CSA $\Lambda$,
un élément $x$ de $\mathbb{R}_{\Lambda}^{n,m}$ est noté
$x=(y,\theta)=(y_{1},y_{2},...,y_{n},\theta_{1},...,\theta_{m})\in\Lambda_{0}^{n}\times\Lambda_{1}^{m}$
avec $y_{i}=\sum_{j=0}^{p}\,y_{i}^{j}e_{j}$ et
$\theta_{k}=\sum_{l=1}^{q}\,\theta_{k}^{l}\eps_{l}.$\\
Nous pouvons définir l'opérateur de Cauchy-Riemann $d''$ dans
$\mathbb{R}_{\Lambda}^{n,m}$ en posant $d''=d_{y}''+d_{\theta}''$,
où $d''_y$ et $d''_\theta$ sont définis en (\ref{d''lambda0^n}) et
(\ref{d"lambda1^m}) :
\begin{equation}\label{d''global}
d''=\sum_{i=1}^{n}\sum_{j=1}^{p}\,dy_{i}^{j}\left(\frac{\partial}{\partial
y_{i}^{j}}-e_{j} \frac{\partial}{\partial y_{i}^{0}}\right)\,
+\,\sum_{l=1}^{m}\sum_{k=0}^{r}\sum_{t=s_{k}+1}^{s_{k+1}-1}\,d\theta_{l}^{t}
(\frac{\partial}{\partial\theta_{l}^{t}}-\frac{\partial}{\partial\theta_{l}^{s_{k}}}a_{t}).
\end{equation}

\section{Solution fondamentale de $d''$ et formule de représentation
intégrale.}

Si $D$ est un ouvert de $\mathbb{R}_{\Lambda}^{n,m}$ borné et à
frontière lisse, si $f$ est une forme différentielle continue dans
$\overline{D}$ et de classe $C^{1}$ dans $D$, nous cherchons une
représentation intégrale de $f$, c'est-à-dire nous essayons
d'écrire
\begin{eqnarray*}
f(x)= \int_{\partial D}\,f(u)K(u,x)\,-\,\int_{D}\,d''f(u)\,. K(u,x)\,+\,d''\int_{D}\,f(u)K(u,x)\,,\;x\in D
\end{eqnarray*}

Soulignons que les formes différentielles intervenant, à
coefficients à valeurs dans $\Lambda$, sont des formes en les
coordonnées {\it réelles} $y_i^j , w_i^j$, ${i= 1,..., n}$,  $j=0
,..., p$\,,\,$\theta_i ^\ell , \tau_i^\ell$, $i=1 ,... ,m$, $\ell
=1,..., q$  de $x=(y,\theta),u=(w,\tau) \in {\Lambda_0}^n \times
{\Lambda_1}^m \approx \R^N$ où $N=n(p+1)+mq$, et les intégrales
sont ainsi définies sans ambiguïté.

\smallbreak\noindent Une telle représentation intégrale peut
s'obtenir dès lors  que $K$ est à coefficients localement
intégrables sur $\mathbb{R}_{\Lambda}^{n,m} \times
\mathbb{R}_{\Lambda}^{n,m}$ et que nous avons en termes de
courants $d''K(u,x)=[\Delta]$ où $[\Delta]$ est le courant
d'intégration sur la diagonale $\Delta=\{(x,x):\,x\in
\mathbb{R}_{\Lambda}^{n,m}\}$ (voir [HP] par exemple). Un noyau
$K(u,x)$ vérifiant cette dernière égalité est appelé noyau
fondamental de l'opérateur de Cauchy-Riemann $d''$. \\
Plus
précisément, soit
$\Psi:\,\overline{D}\times\overline{D}\rightarrow\mathbb{R}_{\Lambda}^{n,m}$
définie par $\Psi(u,x)=u-x$. Alors $[\Delta]=\Psi^{\ast}([0])$ où
$[0]$ est le courant d'évaluation en
$0\in\mathbb{R}_{\Lambda}^{n,m}$. Donc $d''_{u,x} K
(u,x)=[\Delta]$ s'écrit aussi $d''_x \Omega(x) =[0]$ si l'on pose
$K=\Psi^{\ast}\Omega$ ;  chercher une formule de représentation
intégrale pour $d''$ appelle donc à rechercher une solution
fondamentale pour l'opérateur $d''$.

\subsection{Solution fondamentale de $d''$} Comme dans le paragraphe précédent, nous allons d'abord
supposer que f est une forme en $y\in\Lambda_{0}$ (donc $n=1$ dans
un premier temps) et fixer les variables
$\theta\in\Lambda_{1}^{m}$. Dans le cas $\Lambda_{0}=\mathbb{C}$,
une solution fondamentale de $d''\Omega=[0]$ est
$\Omega=\frac{dz}{2i\pi z}$. Mais pour $p>1$, notre algèbre
commutative  $\Lambda_0$ n'est pas intègre (cf. \cite{BM} ou le
théorème de Frobenius classifiant les algèbres réelles
associatives de division de dimensions finies). Nous n'avons pas
non plus de conjugaison sur $\Lambda_0$ comme en analyse complexe
ou quaternionique, et un
noyau résolvant de type Bochner-Martinelli n'est pas envisageable.\\

Nous devons chercher
une solution fondamentale sous une autre forme. \\
Nous munissons $\Lambda_{0}$ d'une norme euclidienne $\|\;\|$ pour
laquelle $(e_{0},...,e_{p})$ est une base orthonormée (avec
$e_{0}$ l'élément unité).\\
Nous sous-entendrons souvent les symbôles $\wedge$ dans les formes
différentielles et utiliserons les notations :

\begin{equation}\label{chapeaux-lambda0}
\widehat{dy{i}}= dy^0 \dots dy^{i-1}dy^{i+1}\dots dy^p \,,\mbox{ et pour  }
i<j\, : \,\widehat{dy^{i}dy^j}= dy^0 \dots
dy^{i-1}dy^{i+1}\dots dy^{j-1}dy^{j+1}\dots dy^p
\end{equation}

\begin{lemma}\label{lemme-sol-fonda}
 Soit $\Omega(y)=\frac{A(y)}{\|y\|^{p+k}}$ où
$A(y)=\sum_{i=0}^{p}\,(-1)^{i}A_{i}(y)\widehat{dy^{i}}$ est une
p-forme à coefficients $A_{i}(y)$ à valeurs dans $\Lambda$ et
polynômes (en les variables réelles $y^{0},...,y^{p}$) homogènes
de degré k et vérifiant $d''\Omega=0$ hors de 0. Alors
$$d''\Omega=\left(\int_{\|y\|\leq1}\,d''A\right)\times[0].$$
\end{lemma}

\noindent  {\sc Preuve :} Soit $\varphi\in
 C_{0}^{\infty}(\Lambda_{0})$ une fonction test. Remarquons que
 pour tout $i=1,\dots ,p$ : $dy^i \wedge \widehat{dy^0}=0$ ; nous
 avons
\begin{eqnarray*}
\langle  d''\Omega,\,\varphi\rangle &=&(-1)^{p+1}\langle\Omega,\,d''\varphi\rangle
=(-1)^{p+1}\langle\Omega,\,\sum_{i=1}^{p}\,
\left(\frac{\partial \varphi}{\partial y^{i}}-e_{i}\frac{\partial \varphi}{\partial y^{0}}\right)dy^{i}\rangle \\
&=& -\int\sum_{i=1}^{p}\,\frac{A_{i}}{\|y\|^{p+k}}\left(\frac{\partial \varphi}{\partial y^{i}}-e_{i}
\frac{\partial \varphi}{\partial y^{0}}\right)dy^{0}...dy^{p} \\
&=& -\lim_{\eps\rightarrow0}\int_{\|y\|\geq\eps}\,\sum_{i=1}^{p}\,\frac{A_{i}}{\|y\|^{p+k}}\left(\frac{\partial \varphi}
{\partial y^{i}}-e_{i}\frac{\partial \varphi}{\partial y^{0}}\right)dy^{0}...dy^{p}
\end{eqnarray*}
\begin{eqnarray*}
&=& -\lim_{\eps\rightarrow0}\int_{\|y\|>\eps}\,d\left[\Big(\sum_{i=1}^{p}\,\frac{(-1)^{i}A_{i}
\widehat{dy^{i}}}{\|y\|^{p+k}}-\frac{A_{i}e_{i}\widehat{dy^{0}}}{\|y\|^{p+k}}\Big)\varphi\right]\,-\,
   ( d''\Omega ) \varphi \\
& =&  \lim_{\eps\rightarrow0}\frac{1}{\eps^{p+k}}\int_{\|y\|=\eps}\,\Big(\sum_{i=1}^{p}A_{i}
\widehat{dy^{i}}(-1)^{i}-A_{i}e_{i}\widehat{dy^{0}}\Big)\varphi
\end{eqnarray*}
\begin{eqnarray*}
&=& \Big(\int_{\|y\|=1}\,\sum_{i=1}^{p}\,(-1)^{i}A_{i}\widehat{dy^{i}}-e_{i}A_{i}\widehat{dy^{0}}\Big)\varphi(0) \\
&=&\int_{\|y\|\leq1}\,\left(\sum_{i=1}^{p}\,\frac{\partial A_{i}}{\partial y^{i}}-e_{i}
\frac{\partial A_{i}}{\partial y^{0}}\right)dy^{0}...dy^{p}.\varphi(0)
=\Big(\int_{\|y\|\leq1}\,d''A\Big)\varphi(0).
\end{eqnarray*}
\noindent en utilisant  $d'' A = d'' (\sum_{i=1}^p \, A_i
\widehat{dy^i})$ pour la dernière égalité. \qed

\medbreak  Ainsi, grâce à ce lemme, pour obtenir une solution
fondamentale pour $d''$, nous allons chercher
$\Omega(y)=\frac{A(y)}{\|y\|^{p+1}}$ avec $A(y)$ p-forme à
coefficients à valeurs dans $\Lambda$ et polynômes homogènes de
degré 1 telle que
$\int_{\|y\|\leq 1}\,d''A\,$ soit inversible.\\

Nous  considérons
\begin{eqnarray*}
\Omega=\frac{1}{\|y\|^{p+1}}\sum_{j=0}^{p}\,(-1)^{j}\widehat{dy^{j}}\sum_{k=0}^{p}\,b_{j}^{k}y^{k}\,=\,\frac{A}{\|y\|^{p+1}}.
\end{eqnarray*}
où $b_{j}^{k}\in\Lambda$ et $x=\sum_{k=0}^{p}\,y^{k}e_{k}.$ Nous
allons choisir les éléments $b_{j}^{k}$ de $\Lambda$  de façon à
obtenir $d''\Omega=0$ hors de 0.\\
 Nous pourrons choisir librement
les $b_0^k$ puisque $d'' \big(\phi(y) \widehat{dy^{0}}\big)=0 $
sur un ouvert  $U \subset \Lambda_0$ dès que $\phi : U\subset
\Lambda_0\rightarrow \Lambda$ est $\R$-différentiable sur $U$.\\
Par un calcul direct :
\begin{eqnarray*}
d''\Omega=\left[\sum_{j=1}^{p}\,\frac{b_{j}^{j}-e_{j}b_{j}^{0}}{\|y\|^{p+1}}\,-\,(p+1)\frac{y^{j}e_{0}-y^{0}e_{j}}{\|y\|^{p+3}}
\sum_{k=0}^{p}\,b_{j}^{k}y^{k}\right]dy^{0}dy^{1}...dy^{p}.
\end{eqnarray*}
$d''\Omega=0$ s'écrit donc
\begin{eqnarray*}
\sum_{j=1}^{p}\,(b_{j}^{j}-e_{j}b_{j}^{0})\,-
\,\frac{p+1}{\|y\|^{2}}\left[\sum_{j=1}^{p}\sum_{k=0}^{p}\,y^{j}y^{k}b_{j}^{k}-\sum_{j=1}^{p}
\sum_{k=0}^{p}y^{0}y^{k}e_{j}b_{j}^{k}\right]=0
\end{eqnarray*}
c'est-à-dire
\begin{eqnarray*}
\scriptstyle{\sum_{j=1}^{p}\,(b_{j}^{j}-e_{j}b_{j}^{0})\,-\,\frac{p+1}{\|y\|^{2}}\left[-({y^{0}})^{\;2}\sum_{j=1}^{p}
\,e_{j}b_{j}^{0}+\sum_{j=1}^{p}\,({y^{j}})^{\;2}b_{j}^{j}+\sum_{j=1}^{p}\,y^{0}y^{j}(b_{j}^{0}-
 \sum_{k=1}^{p}\,e_{k}b_{k}^{j})+\sum_{1\leq j<k}\,y^{j}y^{k}(b_{j}^{k}+b_{k}^{j})\right]=0.}
\end{eqnarray*}
Il faut donc annuler le polynôme obtenu en multipliant par
$\|y\|^{2}$ le premier membre, ce qui donne:
\begin{equation}\label{lesbjk}
\begin{cases}
\; b_{j}^{0}=\sum_{k=1}^{p}\,e_{k}b_{k}^{j};\;j=1,2,...,p \\
\; b_{j}^{k}+b_{k}^{j}=0;\;j,k=1,...,p;\,j\neq  k \\
\;
b_{j}^{j}=-\sum_{k=1}^{p}\,e_{k}b_{k}^{0}=\frac{1}{p+1}\sum_{k=1}^{p}\,(b_{k}^{k}-e_{k}b_{k}^{0});\;j=1,...,p.
\end{cases}
\end{equation}
Notons $b\,=\,-\sum^{p}_{k=1}\,e_{k}b^{0}_{k}$ la valeur commune
des $b^{j}_{j}$, $j=1,\dots, p$. La première ligne de
(\ref{lesbjk}) permettra de définir $b_{j}^{0}$ quand nous
connaîtrons les $b_{k}^{j}$. La deuxième ligne donnera $b_{k}^{j}$
quand nous connaîtrons $b_{j}^{k}$ pour $j<k$. Compte tenu de ces
deux lignes :
\begin{eqnarray*}
b&=& -\sum_{k=1}^{p}\,e_{k}\sum_{l=1}^{p}\,e_{l}b_{l}^{k}
= -\sum_{k,l=1}^{p}\,e_{k}e_{l}b_{l}^{k}
 \,=\,-\sum_{k=1}^{p}\,e_{k}^{\;2}b_{k}^{k}\,-\,\sum_{k<l}\,e_{k}e_{l}b_{l}^{k}\,
 -\,\sum_{l<k}\,e_{k}e_{l}b_{l}^{k} \\
 &=&-\sum_{k=1}^{p}\,e_{k}^{\;2}b\,-\,\sum_{k<l}\,e_{k}e_{l}(b_{k}^{l}+b_{l}^{k})
 \, =\,-b\sum_{k=1}^{p}\,e_{k}^{\;2}.
\end{eqnarray*}
ou encore
\begin{eqnarray}
b\sum_{k=0}^{p}\,e_{k}^{\;2}=0.
\end{eqnarray}
Alors trois cas se présentent. \\
Si $\sum_{k=0}^{p}\,e_{k}^{\;2}$ n'a pas de diviseur de 0,
nécessairement $b=0$ et donc $\int_{\|x\|\leq1}\,d''A\,=\,0$: nous
n'avons pas de représentation intégrale. \\
Si $\sum_{k=0}^{p}\,e_{k}^{\;2}$ non nul a des diviseurs de 0, on
prend pour b l'un de ces diviseurs et alors
$\int_{\|x\|\leq1}\,d''A$ est non inversible: nous n'avons pas une
bonne représentation intégrale. \\
Enfin, si $\Lambda_{0}$ vérifie la condition ($A_{0}$) suivante :
\begin{equation}\label{A0}
\quad\bf{(A_{0})} \;   \mbox{ il existe une base
$(e_{0}=1,e_{1},...,e_{p})$ de $\Lambda_{0}$ vérifiant
$\displaystyle{\sum_{k=0}^{p}\,e_{k}^{\;2}=0}$ \,,}
\end{equation}

\noindent nous prenons pour b un élément inversible de $\Lambda$,
par exemple $b=e_{0}$, et choisissons $b_{j}^{k}=0$ si $j\neq k$,
$1\leq j\leq p$ ; nous déduisons alors du lemme
\ref{lemme-sol-fonda}, puisque
 $\int_{\|x\|\leq 1}\,d''A\,=\,(p+1)e_{0}Vol(B(0,1))$ où
Vol(B(0,1)) est le volume de la boule unité de $\Lambda_{0}$ :
$$\Omega=\frac{1}{(p+1)Vol(B(0,1))\|y\|^{p+1}}\big( \sum_{j=1}^{p}\,(-1)^{j}(y^{0}e_{j}+y^{j}e_{0})\widehat{dy^{j}} +
\sum_{k=0}^p b_0^k y^k \widehat{dy^{0}}\big)$$
est une solution fondamentale de $d''$ dans $\Lambda_0$.\\
Nous allons choisir les $b_0^k$ en sorte que $d' \Omega =0$, ce
qui simplifiera d'une part certains calculs d'autre part la
formule de représentation intégrale  ; ceci conduit à $b_0^0 =
e_0$ et $b_0^k = - e_k$ si $k\geq 1$. La forme $dY=\sum_{k=0}^p
e_k dy^k$ est en facteur dans $A$ et nous obtenons, avec les
notations (\ref{chapeaux-lambda0}) :
\begin{proposition}\label{Omega0}
$$\Omega_0=\frac{1}{(p+1)Vol(B(0,1))\|y\|^{p+1}}\sum_{j=1}^{p}\,(-1)^{j}(y^{0}e_{j}+y^{j}e_{0})
(\sum_{k=0}^p e_k dy^k )\widehat{dy^0 dy^{j}}$$
est une solution fondamentale de $d''$ dans $\Lambda_0$.
\end{proposition}

\smallbreak\noindent {\sc Exemples }:  La condition ($A_0$) est
vérifiée dans les exemples \ref{exemple-C} et \ref{exemple-table}.
%et, lorsque $k
%\equiv 2$ modulo 4, dans l'exemple 3 du paragraphe 1.
 Elle ne
l'est pas dans le second exemple (ce serait d'ailleurs en
contradiction avec la proposition
\ref{harmonicite} du paragraphe 4).\\

 {\it La condition $( A_{0} )$ sera toujours supposée vérifiée
dans la suite}.

\medbreak

$\bullet$ cherchons maintenant une solution fondamentale pour
l'opérateur $d''$ défini sur $\Lambda_{1}$.
 De façon à disposer d'une expression explicite pour $d''$, nous supposerons satisfaite la condition (\ref{préA1}) ;
 par conséquent, nous avons :
$$d''=\sum_{k=1}^{r}\sum_{j=s_{k}+1}^{s_{k+1}-1}d\theta^{j}\left(\frac{\partial}{\partial\theta^{j}}-a_{j}\frac{\partial}{\partial\theta^{s_{k}}}\right),$$
Nous munissons $\Lambda_{1}$ d'un produit scalaire en supposant
que $(\eps_{1},...,\eps_{q})$ est une base orthonormée, notons
$\|\;\|$ la norme associée, et cherchons une solution fondamentale
de la forme:
\begin{eqnarray*}
\Omega=\sum_{k=1}^{r}\sum_{j=s_{k}}^{s_{k+1}-1}\,(-1)^{j-1}\widehat{d\theta^{j}}
\sum_{l=1}^{q}\frac{b_{j}^{l}\theta^{l}}{\|\theta\|^{q}}
\end{eqnarray*}
On calcule $d''\Omega$ hors de 0, et l'on note
$dV=d\theta^{1}\wedge ...\wedge d\theta^{q}$ ; comme $d''
\big(\phi(y,\theta) \widehat{d\theta^{s_k}}\big)=0$ sur tout
ouvert $U$ sur lequel $\phi$ est $\R$-différentiable, nous
pourrons choisir arbitrairement les $b^l_{s_k}$.
\begin{eqnarray*}
d''\Omega\;=\;\frac{1}{\|\theta\|^{q}}\sum_{k=1}^{r}\sum_{j=s_{k}+1}^{s_{k+1}-1}\left(b_{j}^{j}-a_{j}b_{j}^{s_{k}}\right)dV
\;-\;\frac{q}{\|\theta\|^{q+2}}\sum_{k=1}^{r}\sum_{j=s_{k}+1}^{s_{k+1}-1}\left(\theta^{j}-a_{j}\theta^{s_{k}}\right)
\left(\sum_{l=1}^{q}\,b_{j}^{l}\theta^{l}\right)dV.
\end{eqnarray*}
donc
\begin{eqnarray*}
\|\theta\|^{q+2}d''\Omega \;{\widehat{dV}}\;&=&\;\|\theta\|^{2}\sum_{k=1}^{r}\sum_{j=s_{k}+1}^{s_{k+1}-1}
\left(b_{j}^{j}-a_{j}b_{j}^{s_{k}}
\right)\;\\
&& -\,q \left(-\sum_{k=1}^{r}(\theta^{s_{k}})^{2}\sum^{s_{k+1}-1}_{j=s_{k}+1}a_{j}b^{s_{k}}_{j}\,+\,\sum^{r}_{k=1}
\sum^{s_{k+1}-1}_{j=s_{k}+1}\,(\theta^{j})^{2}b^{j}_{j}\right)\\
&& + q \sum_{s_{k}< s_{k'}}\,\theta^{s_{k}}\theta^{s_{k'}}
\left(\sum^{s_{k+1}-1}_{j=s_{k}+1}\,a_{j}b^{s_{k'}}_{j}\,+\,\sum^{s_{k'+1}-1}_{j=s_{k'}+1}\,a_{j}b^{s_{k}}_{j}\right)
\;\\
&& - q\Big[\sum^{r}_{k=1}\sum^{r}_{k'=1}\sum^{s_{k'+1}-1}_{j=s_{k'}+1}\,\theta^{s_{k}}\theta^{j}\left(b^{s_{k}}_{j}\,-\,
\sum^{s_{k+1}-1}_{l=s_{k}+1}\,a_{l}b^{j}_{l}\right)\,+ \sum_{\substack{ j<k\\j,k\neq s_{1},...,s_{r}}}
\,\theta^{j}\theta^{k}(b^{k}_{j}\,+\,b^{j}_{k})\Big].
\end{eqnarray*}
La nullité de  $d''\Omega$ hors de 0  s'écrit donc:
\begin{equation*}
\begin{cases}
& \,b^{k}_{j}+b^{j}_{k}\,= \,0\;si\;j\neq k,\;j\,et\,k\notin\{s_{1},...,s_{r}\}\\
&\, b^{s_{k}}_{j}\,=\,\displaystyle{\sum^{s_{k+1}-1}_{l=s_{k}+1}}\,a_{l}b^{j}_{l}\\
&\, 0 = \displaystyle{\sum^{s_{k+1}-1}_{j=s_{k}+1}\,a_{j}b^{s_{k'}}_{j}\,+\,\sum^{s_{k'+1}-1}_{j=s_{k'}+1}\,a_{j}b^{s_{k}}_{j}}\quad pour \; k\neq k'\\
&\, \displaystyle{\sum^{r}_{k=1}\sum^{s_{k+1}-1}_{j=s_{k}+1}\,(b^{j}_{j}\,-\,a_{j}b^{s_{k}}_{j})\,=\,qb^{j}_{j}\,=\,
-q\sum^{s_{k_{0}+1}-1}_{j=s_{k_{0}}+1}\,a_{j}b^{s_{k_{0}}}_{j}}\mbox{ pour tout } \;k_{0}\in \{1,...,r\}.
\end{cases}
\end{equation*}

La troisième ligne est une conséquence des deux premières ; en
effet:
\begin{eqnarray*}
\sum^{s_{k+1}-1}_{j=s_{k}+1}\,a_{j}b^{s_{k'}}_{j}\,+\,\sum^{s_{k'+1}-1}_{j=s_{k'}+1}\,a_{j}b^{s_{k}}_{j}\,
&=&\,\sum^{s_{k+1}-1}_{j=s_{k}+1}\sum^{s_{k'+1}-1}_{l=s_{k'}+1}\,a_{j}a_{l}b^{j}_{l}\,
+\,\sum^{s_{k'+1}-1}_{j=s_{k'}+1}\sum^{s_{k+1}-1}_{l=s_{k}+1}\,a_{j}a_{l}b^{j}_{l}\\
&=&\,\sum^{s_{k+1}-1}_{j=s_{k}+1}\sum^{s_{k'+1}-1}_{l=s_{k'}+1}\,a_{j}a_{l}(b^{j}_{l}\,+\,b^{l}_{j})\,=\,0.
\end{eqnarray*}
D'après la quatrième ligne tous les
$b^{j}_{j}$ sont égaux à une constante $b\in\Lambda$. Dès lors, la
dernière égalité de celle-ci s'écrit, en tenant compte de la seconde ligne:
$$b\,=\,-\sum^{s_{\scriptscriptstyle{k+1}}-1}_{j,l=s_{k}+1} \,a_{j}a_{l}b^{j}_{l}\\
=\,-\sum_{j<l}\,a_{j}a_{l}(b^{l}_{j}+b^{j}_{l})\,-\,\sum^{s_{k+1}-1}_{j=s_{k}+1}\,a^{\;2}_{j}b,\;\forall k=1,\dots,r$$
c'est-à-dire:
\begin{equation}\label{égalitésA1}
b\sum^{s_{k+1}-1}_{j=s_{k}}\,a^{\;2}_{j}\,=\,0 \qquad \forall k=1,\dots,r
\end{equation}
La première égalité de la quatrième ligne nous conduit de même à
(\ref{égalitésA1}).

On veut, tout comme dans $\Lambda_{0}$, que b soit un élément
inversible de $\Lambda$, afin d'obtenir une bonne formule de
représentation intégrale ; nous sommes ainsi amenés à imposer à
$\Lambda_{1}$, en tenant compte de (\ref{préA1}) la condition
suivante :

\begin{equation}\label{A1}
\bf{(A_{1})}\begin{array}{l}
\mbox{ il existe une base $(\eps_{1},...,\eps_{q})$ de
$\Lambda_{1}$ et une suite finie
$s_{1}= 1<s_{2}<...<s_{r}<s_{r+1}=q+1$} \\
\mbox{  telles que, pour tout
$j=1,...,q$,
  il existe $a_{j}\in\Lambda_{0}$ vérifiant
$\eps_{j}=a_{j}\eps_{s_{k}}\,$ si $\,s_{k}\leq j<s_{k+1}$,}\\
\mbox{  avec
$a_{s_{1}}=a_{s_{2}}=...=a_{s_{r}}=e_{0}$\;
et\;
$\sum_{j=s_{k}}^{s_{k+1}-1}\,a_{j}^{\;2}=0\;$ pour tout $k=1,...,r$}\,.
\end{array}
\end{equation}

Dès lors que la condition ($A_{1}$) est satisfaite, nous pouvons
écrire une solution fondamentale explicite pour l'opérateur $d''$
sur $\Lambda_1$ ; par exemple, prenant $b=e_{0}$ et $b_{j}^{k}=0$
si $j\neq k\,,\;j,k\geq2$, nous obtenons comme candidat :
\begin{equation}\label{Omega-non-fermé}
\Omega\,=\,\frac{1}{q\,Vol(B(0,1)\, \|\theta\|^{q}}\sum^{r}_{k=1}\big(\sum^{s_{k+1}-1}_{j=s_{k}+1}\,(-1)^{j-1}
(\theta^{j}e_{0}+\theta^{s_{k}}a_{j})\widehat{d\theta^{j}}\,  + (\sum_{l=1}^q b^l_{s_k}\widehat{d\theta^{s_k}}\big)
\end{equation}
Nous imposons de plus, comme dans l'étude concernant $\Lambda_0$,
$d' \Omega =0$, ce qui conduit à choisir pour tout $k=1,\dots ,r$
: $ b_{s_k}^{s_k}=e_0 \,;\, b^{l}_{s_k}= - a_l \mbox{ si } s_k <
l< s_{k+1}\,;\,b^{l}_{s_k}= 0 \mbox{ si } l < s_k \mbox{ ou }
s_{k+1}< l$.\\
Nous obtenons alors, avec ces choix, le noyau suivant :
\begin{equation}\label{Omega1}
\Omega_1\,=\,\frac{1}{q\,Vol(B(0,1)\, \|\theta\|^{q}}\sum^{r}_{k=1}\sum^{s_{k+1}-1}_{j=s_{k}+1}\,(-1)^{j-1}
(\theta^{j}e_{0}+\theta^{s_{k}}a_{j})\,(\sum_{\ell =s_k }^{s_{k+1}-1} a_\ell d\theta^\ell )\,\widehat{d\theta^{s_k} d\theta^{j}}\,
\end{equation}
 Pour
vérifier qu'il s'agit bien d'une solution fondamentale de $d''$
dans $\Lambda_1$, nous devons étudier la singularité en 0, ce qui
se fait de façon identique au cas de $\Lambda_{0}$, à l'aide de
l'analogue du lemme \ref{lemme-sol-fonda} (voir aussi, ci-dessous,
la démonstration dans $\Lambda_{0}^{n}\times\Lambda_{1}^{m}$, dont
le cas présent est un cas particulier).

\smallbreak\noindent {\sc Exemple.} Si l'on prend
$\Lambda_{1}=\eps_{1}\Lambda_{0}$, alors la condition ($A_{1}$) se
ramène à la condition
 ($A_{0}$) si l'on pose $\eps_{j}=\eps_{1}e_{j-1}$. \\

\medbreak {\it La condition ($A_{1}$) sera supposée satisfaite
dans toute la suite}.

\medbreak \noindent Désormais, il reste à chercher une solution
fondamentale pour $d''$ défini sur un superespace
$\mathbb{R}_{\Lambda}^{n,m}=\Lambda_{0}^{n}\times\Lambda_{1}^{m}$.
Nous supposons toujours satisfaites les conditions ($A_{0}$) et
($A_{1}$) en sorte que l'opérateur de Cauchy-Riemann sur
$\mathbb{R}_{\Lambda}^{n,m}$ est défini comme en
(\ref{d''global}).

\medbreak Nous allons utiliser les résultats précédents obtenus pour
$\Lambda_0$ et $\Lambda_1$. Notons

\begin{equation}\label{formes-eta-omega}
\begin{cases}
\omega (y_i)=& dy_i^0\dots dy_i^p \, (= dy_{i}^{0}\wedge \dots \wedge dy_{i}^{p}) \;\qquad \mbox{ pour }\; i=1,\dots, n\\
\lambda (\theta_i ) =& d\theta_{i}^{1} \dots  d\theta_{i}^{q} \qquad \qquad\qquad\qquad\qquad  \mbox{ pour }\; i=1,\dots, m\\
\eta (y_i )=& \displaystyle{\sum_{j=1}^{p}\,(-1)^{j}(y_{i}^0 e_{j}+y_{i}^j e_{0})dY_i \, dy_i^1 \dots dy_i^{j-1} dy_i^{j+1}\dots dy_i^p
\qquad\mbox{ pour }\; i=1,\dots, n}\\
& \qquad\qquad\quad\qquad\quad\qquad\quad\mbox{ où } dY_i = \sum_{k=0}^p e_k
dy_i^k \\
\nu(\theta_i )= &\displaystyle{\sum^{r}_{k=1}\sum^{s_{k+1}-1}_{j=s_{k}+1}\,(-1)^{j-1}
(\theta_i^{j} e_{0}+\theta_i^{s_{k}}a_{j})\,dZ_k (\theta_i ) \,d\theta_i^1 \dots d\theta_i^{s_k -1}d\theta_i^{s_k +1}
\dots d\theta_i^{j-1}
d\theta_i^{j+1} \dots d\theta_i^{q}}  \\
&\quad\qquad\qquad\quad\qquad \mbox{ pour } \;  i=1,\dots, m\;,\; \mbox{ où }
  dZ_k (\theta_i ) =\sum_{\ell =s_k }^{s_{{k+1}} -1} a_\ell d\theta_i^\ell \;,\; k=1,\dots ,r \\
dV \;=&  \displaystyle{\bigwedge_{j= 1}^n
\omega (y_i )  \bigwedge_{\ell=1}^m \lambda
(\theta _\ell )}
\end{cases}
\end{equation}

\medbreak\noindent Comme les formes différentielles $dY_i$ et
$dZ_k (\theta_i)$ définies dans (\ref{formes-eta-omega}) sont à
coefficients dans $\Lambda_0$, nous avons :

\qquad $dY_i \wedge dY_i =0$, $\forall i=1,\dots ,n$ et $dZ_k (\theta_i) \wedge dZ_k (\theta_i)=0$, \;$\forall i=1,\dots,m\,,\, \forall k=1,\dots , r$.

\medbreak\noindent Tenant compte de la condition ($A_0$) (cf
(\ref{A0})), nous pouvons écrire pour tout $i=1,\dots ,n$ :
\begin{equation}\label{eta}
\begin{array}{l}
d \eta (y_i ) = d'' \eta (y_i ) = (p+1) e_0 \,\omega (y_i )\\
 \frac{1}{2}{\displaystyle{d ||x||^{2} \wedge \eta (y_i
)\bigwedge_{{\substack{{j =1}\\{j\neq i}}}}^n
 \omega (y_j )}  \bigwedge_{\ell=1}^m \lambda
(\theta _\ell)}
 = \frac{1}{2}
{\displaystyle{ d''_{y_i} (||x||^{2}) \wedge \eta (y_i )\bigwedge ^n_{{\substack{{j =1}\\{j\neq i}}}}
 \omega (y_j )
 \bigwedge_{\ell=1}^m \lambda (\theta _\ell)}}\\
 = (-1)^{(p+1)^2 (i-1)}
 ||y_i||^2  e_0 \,dV
 \end{array}
\end{equation}

\noindent Tenant compte de ($A_1$) ((cf (\ref{A1})), nous obtenons
pour tout $i=1, \dots , m$ :
\begin{equation}\label{nu}
\begin{array}{l}
d \nu (\theta_i ) = d'' \nu (\theta_i ) = q e_0 \, \lambda (\theta_i )\mbox{ et }\\
  \frac{1}{2}{\displaystyle{d ||x||^{2} \wedge \nu(\theta_i )\bigwedge_{j= 1}^n
\omega (y_i )
 \bigwedge_{\ell\neq i} \lambda
(\theta _\ell )}}
 =\frac{1}{2}
\displaystyle{ d''_{\theta_i} (||x||^{2})\wedge \nu(\theta_i ) \bigwedge_{j=
1}^n \omega (y_i )
 \bigwedge_{\ell\neq i} \lambda
(\theta _\ell)} \\
 = (-1)^{nq(p+1)+q^2 (i-1)}\,
 ||\theta_i||^{2}  e_0 \, dV
 \end{array}
\end{equation}

Posons :
\begin{eqnarray}\label{def-Atilde}
 &&\Tilde{\Omega}=\frac{\Tilde{A}}{\|x\|^{n(p+1)+qm}} \qquad \mbox{ où }\nonumber \\
 &&\begin{array}{c}
 \Tilde{A}=\displaystyle{\sum_{i=1}^{n}(-1)^{\scriptscriptstyle{(p+1)^2 (i-1)}}\,\eta (y_i
)\bigwedge_{\substack{ j =1 \\j\neq i}}^n
 \omega (y_j )  \bigwedge_{\ell=1}^m \lambda
(\theta _\ell )
+\,\sum_{i=1}^{m} (-1)^{\scriptscriptstyle{nq(p+1)+q^2 (i-1)}}\nu(\theta_i )\bigwedge_{j= 1}^n
\omega (y_j )
 \bigwedge_{\substack{\ell =1\\\ell\neq i}}^m \lambda
(\theta_\ell )}\\
= \displaystyle{\sum_{i=1}^{n} \,dY_i
 \sum_{j=1}^p (-1)^j (y_i^j e_0 +y_i^0 e_j){\widehat{dy_i^0 dy_i^j}} \,+ \,\sum_{i=1}^{m}\sum_{k=1}^{r} dZ_k (\theta_i )
\sum^{s_{k+1}-1}_{\ell=s_{k}+1 } (-1)^{\ell -1}(\theta_{i}^{\ell} e_{0}
+ \theta_{i}^{s_k} a_{\ell})\widehat{d\theta_{i}^{s_k}d\theta_i^\ell}}
\end{array}
\end{eqnarray}

 \noindent où l'on a noté,  en accord avec (\ref{formes-eta-omega}) :
\begin{equation}\label{nota-chapeaux}
\begin{array}{l}
 \widehat{dy_{i}^{0}dy_{i}^{j}}= \;\omega (y_{1})\dots\omega (y_{i-1})dy_{i}^{1}\dots dy_{i}^{j-1}dy_{i}^{j+1}...dy_i^p
 \omega(y_{i+1})\dots \omega (y_n )\lambda (\theta_1 )\dots \lambda (\theta_m )\\
\widehat{d\theta_{i}^{s_k}d\theta_{i}^{\ell}}= \omega (y_{1})...\omega (y_{n})\lambda (\theta_1 )\dots  \lambda (\theta_{i-1} )
d\theta_i^1 .. d\theta_i^{s_k -1}d\theta_i^{s_k +1}
.. d\theta_i^{\ell -1}d\theta_i^{\ell +1}...d\theta_i ^q \lambda (\theta_{i+1} )\dots \lambda (\theta_m )
 \end{array}
\end{equation}

\smallbreak\noindent Nous déduisons des calculs (\ref{eta}) et
(\ref{nu}) précédents
\begin{equation}\label{Atilde}
  d'' \tilde{A}=d \tilde{A} = (n(p+1) +mq) e_0 \, dV
  \end{equation}
  \begin{equation}\label{Omega-fermé}
d{\Tilde{\Omega}} =
d''\Tilde{\Omega}= 0 \qquad  \;\mbox{ hors de } \{0\}
\end{equation}
Il reste maintenant à étudier la singularité de $\Tilde{\Omega}$ ;
remarquons que les coefficients de $\Omega$ sont de classe
$C^\infty$ hors de $0$ et localement intégrables  ; si
$\varphi=\varphi(y,\theta)$ est une fonction $C^{\infty}$ à
support compact sur $\mathbb{R}_{\Lambda}^{n,m}$, alors, en notant
$N=n(p+1)+mq$ :

\begin{eqnarray*}
\langle  d''\tilde{\Omega},\varphi\rangle
&=& (-1)^{N}\langle\Tilde{\Omega},
d''\varphi\rangle = (-1)^{N}\langle\Tilde{\Omega},
d\varphi\rangle
 \nonumber
=(-1)^{N} \lim_{\eps\rightarrow 0}\int_{\|x\|\geq\eps}\,\Tilde{\Omega} \wedge d\varphi \nonumber\\
&=& - \lim_{\eps\rightarrow 0}\, \int_{\|x\|\geq \eps}\, d \big(\Tilde{\Omega}\varphi) \qquad\qquad\qquad\;
{\mbox{ d'après  (\ref{Omega-fermé})}}\nonumber \\
&=&  \lim_{\eps\rightarrow 0}\, \frac{1}{\eps^{N}}\, \int_{\|x\|=\eps} \, {\tilde{A}} \varphi \nonumber \\
&= & \lim_{\eps\rightarrow 0}\frac{\varphi (0)}{\eps^{N}}\, \int_{\|x\|=\eps}\, \tilde{A} \qquad\qquad \qquad  \mbox{ car } ||\tilde{A}||
\lesssim ||x|| \nonumber\\
&=&  \lim_{\eps\rightarrow 0}\frac{\varphi (0)}{\eps^{N}}\,  \int_{\|x\|\leq \eps}\, d\tilde{A}
\end{eqnarray*}

\noindent et, en tenant compte de (\ref{Atilde}) :
\begin{equation}\label{Omega-fonda}
\langle  d''\tilde{\Omega},\varphi\rangle =  (n(p+1)+mq)  Vol(B(0,1)) \,e_0\, \varphi (0)
\end{equation}

\noindent Par suite, $\Omega := {[(n(p+1)+qm) Vol(B(0,1))]}^{-1}
\Tilde{\Omega}$ \, est une solution fondamentale pour l'opérateur
$d''$ dans $\mathbb{R}_{\Lambda}^{n,m}$

\smallbreak\noindent {\sc remarque :} On vérifie aisément, grâce à
(\ref{eta}), (\ref{nu})  et (\ref{def-Atilde}) l'égalité au sens
des courants $d' \Omega =0$.

\medbreak Nous pouvons résumer les résultats de ce paragraphe dans
le théorème \ref{thm-sol-fonda} suivant.

\noindent Rappelons que les conditions $(A_0 )$ et $(A_1 )$ sont
définies en (\ref{A0}) et (\ref{A1}), l'opérateur $d''$ sur
$\Lambda_{0}^{n}\times\Lambda_{1}^{m}$ en (\ref{d''global}).

 \begin{theoreme}\label{thm-sol-fonda}
  Si $\Lambda_{0}$
vérifie la condition ($A_{0}$) et $\Lambda_{1}$  la condition
($A_{1}$), l'opérateur de Cauchy-Riemann $d''$ dans
$\mathbb{R}_{\Lambda}^{n,m}$ admet une solution fondamentale
$\Omega$ donnée par :
\begin{eqnarray}
\Omega (x)&=&\frac{c\scriptscriptstyle{(n,m,p,q)}}{\|x\|^{N} }
\Big( \displaystyle{\sum_{i=1}^{n} \,dY_i
 \sum_{j=1}^p (-1)^j (y_i^j e_0 +y_i^0 e_j){\widehat{dy_i^0 dy_i^j}}} \, \nonumber \\
 & + & \,\displaystyle{\sum_{i=1}^{m}\sum_{k=1}^{r} dZ_{k} (\theta_i )
\sum^{s_{k+1}-1}_{\ell=s_{k}+1 } (-1)^{\scriptscriptstyle{(\ell -1 )}}(\theta_{i}^{\ell} e_{0}+ \theta_{i}^{s_k} a_{\ell})\widehat{d\theta_{i}^{s_k}d\theta_i^\ell}}  \Big)
\end{eqnarray}
\noindent où $N=n(p+1)+mq\;,\;c(n,m,p,q)= -\big(N\, \Vol
(B(0,1))\big)^{-1}, \,\widehat{dy_i^0 dy_{i}^{j} }$
 et $\widehat{d\theta_{i}^{s_k}d\theta_i^\ell} $ étant définis en
(\ref{nota-chapeaux})
\end{theoreme}

\subsection{Formule de représentation intégrale pour les fonctions (et les formes différentielles)}

Soit D un ouvert de $\mathbb{R}_{\Lambda}^{n,m}$ borné à frontière
lisse, soit
$\Psi:\,\overline{D}\times\overline{D}\rightarrow\mathbb{R}_{\Lambda}^{n,m}$
définie par $\Psi(x,x')=x-x'$. Définissons, avec les notations (\ref{formes-eta-omega}) et (\ref{nota-chapeaux}) :
\begin{equation}\label{defn-K}
\begin{array}{l}
  K(x,x')={\Psi^{\ast}} \Omega \\
=\dfrac{c\scriptscriptstyle{(n,m,p,q)}}{\|x'-x\|^{N}}
\Big( \scriptstyle{\sum_{i=1}^{n} (dY_i - dY'_i ) \sum_{j=1}^{p}\,(-1)^{\scriptscriptstyle{j}}
\big(({y}_{i}^{j} - {y'}_{i}^{j}) e_{0}+({y}_{i}^{0} -{y'}_{i}^{0})e_j
 \big)
\overset{{\makebox{\put(-31,0){\line(6,1){31}}}
\makebox{\put(0.5,5){\line(6,-1){32}}}
}}{{d(y_{i}^{0}-{y'}_{i}^{0}) d({y}_{i}^{j}-{y'}_{i}^{j})}}}  \\
+\scriptstyle{\sum_{i=1}^{m}\sum_{k=1}^{r} dZ_k (\theta_{i}-{\theta'}_i )
\sum^{s_{k+1}-1}_{l=s_{k}+1}\,(-1)^{\scriptscriptstyle{l-1}}
\big(({\theta}_{i}^{l}-{\theta'}_{i}^{l})e_{0}+({\theta}_{i}^{s_{k}}-{\theta'}_{i}^{s_{k}})a_{l}
\big)
\overset{{\makebox{\put(-33.4,0){\line(6,1){34}}}
\makebox{\put(1,5.9){\line(6,-1){33}}}
}}{{d(\theta_{i}^{s_k}-{\theta'}_{i}^{s_k})
d({\theta}_{i}^{l}-{\theta'}_{i}^{l})}}}
\Big)
\end{array}
\end{equation}

Les coefficients de $K$ sont de classe $C^\infty$ hors de la
diagonale $\Delta$ de $ \R^{n,m}_\Lambda \times \R^{n,m}_\Lambda $
et localement intégrables sur $\R^{n,m}_\Lambda \times
\R^{n,m}_\Lambda $.

\noindent  Le noyau $K$ vérifie au sens des courants, si l'on pose
$d=d_x +d_{x'}$ et $d''=d_{x'}''+d_{x}''$:
\begin{equation}\label{dK}
d''K(x,x')=d''\Psi^{\ast}\Omega (x)=\Psi^{\ast}d''\Omega={\Psi^{\ast}}[0]=[\Delta]
\qquad \qquad \mbox{ et }\;  d K(x',x) = \Psi ^\ast d\Omega = [\Delta]
\end{equation}

\smallbreak Précisons
ce que nous entendons par bidegré d'une forme. Toute
forme différentielle en les N variables réelles $y_i^j$, ${1\leq
i\leq n}$, $0\leq j\leq p$ et $\theta_i^\ell$, ${1\leq i\leq m}$,
${1 \leq j\leq q}$ peut s'écrire de
 manière unique
 en fonction des 1-formes (dites de bidegré (1,0)) de l'ensemble\\
 $\mathfrak{F}  =\{dY_i \,, {{i=1,\dots ,n}}\;;\; dZ_k (\theta_i )\;,{{i=1,\dots ,m \;, k=1,\dots
 ,r}}\}$ et des 1-formes
 dites de bidegré (0,1) appartenant à l'ensemble \\
 $\mathfrak{H}= \{dy^j_i \,,\, {\textstyle{i=1,\dots ,n\,,j=1,\dots ,p}}\, \;;\; d\theta^\ell_i \;,{{i=1,\dots
 ,m \,,\, \ell =1,\dots ,q\,,\, \ell \neq s_k\,,\,\forall
 k=1,..., r}}\,\}$.\\
 Une forme $f$ de degré $d$ est dite de bidegré $(a,d-a)$ si dans chacune de ses
 composantes apparaissent  exactement $a$ 1-formes appartenant à   $\mathfrak{F}$.

\begin{corollary}\label{corol-repres} Sous les hypothèses du théorème
précédent  :

(a)   si f est une fonction   de classe $C^{1}$ dans D, continue
sur $\overline{D}$ ainsi que $df$, alors pour tout $x'\in D$, nous
avons :
\begin{eqnarray*}
f(x')= \int_{\partial D}\, K(x,x') f(x)\;-\;\int_{D}\,d''f(x)\,{\scriptstyle{\wedge}} K(x,x')\,.
\end{eqnarray*}

(b) plus généralement, si $f$ est une forme de degré $d$ de
bidegré $(a,d-a)$, de classe $C^{1}$ dans D, continue sur
$\overline{D}$ ainsi que $df$, nous avons pour tout $x'$ de $D$:
\begin{equation*}
f(x')=\int_{\partial D}\,f(x) K_{(a,d-a)}(x,x')\,+\,(-1)^{d-1}\int_{D}\,d''f(x){\scriptstyle{\wedge}}K(x,x')\,
+ \,(-1)^{d} \,d''\int_{D}\,f(x) K(x,x')
\end{equation*}
\noindent où $K_{(a,d-a)}$ désigne la composante de $K$ de bidegré
$(a,d-a)$ en $x'$.
\end{corollary}

Comme nous utiliserons uniquement le résultat du (a) dans la suite
de l'article, nous reportons en appendice la preuve du (b).

\noindent (a) Nous en proposons ici une preuve élémentaire. Il
suffit, vu la régularité des intégrales, de prouver l'égalité au
sens des distributions. Supposons $f$ de classe $C^\infty$ dans
$D$ ; soit $g(x')$ une
N-forme de classe $C^\infty$ à support compact dans $D$, où  $N$ est la dimension de $\R^{n,m}_\Lambda$. \\
 D'après (\ref{dK}), nous avons hors de $\Delta$  : $ d_{x,x'} \big[ f(x) K(x,x')g(x')
\big] =
 df(x){\scriptstyle{\wedge}} K(x,x')g(x')  $

\noindent Or $ K(x,x') g(x')=  L(x,x') g(x')$ avec
\begin{equation}\label{L}
\begin{array}{l}
\quad L(x,x')=c{\scriptscriptstyle{(n,m,p,q)}}{\|x'-x\|^{-N}\,
\tilde{B}} {\mbox{ où, avec les notations (\ref{formes-eta-omega}) et (\ref{nota-chapeaux}) }}\\
\quad \tilde{B} =  \sum_{i=1}^{n} dY_i
\sum_{j=1}^{p}\,
b_{i,j}(y_i ,y'_i )\,  \widehat{dy_i^0 dy_i^j}
\,
+\,\sum_{i=1}^{m}\sum_{k=1}^{r} dZ_k (\theta_i) \sum^{s_{k+1}-1}_{l=s_{k} +1}\,
\beta_{k,l}(\theta_i , \theta'_i )\,
\widehat{ d\theta_i^{s_k}  d\theta^l_i  }
\end{array}
\end{equation}

\noindent les $b_{i,j} $ et les $\beta_{k,l}$ étant des fonctions
à valeurs dans $\Lambda$ (cf. la définition (\ref{defn-K}) de
$K$).

\noindent Pour tout $\ell =1,...,n$ : $dY_\ell \wedge \tilde{B}=0$
et pour tout $\ell =1, \dots , m\,,k=1,...,r$ : $dZ_k (\theta_\ell
)\wedge \tilde{B}=0$, par suite $d'f(x) \wedge L(x,x')    =0$

\noindent Nous avons donc hors de $\Delta$ : $ d_{x,x'} \big[ f(x)
K(x,x')g(x')\big] =
 d''f(x){\scriptstyle{\wedge}} K(x,x')g(x') $.
Par la formule de Stokes, nous obtenons pour $0<\eps <<1$ :
\begin{eqnarray*}
\int_{{\partial D}  \times D}\,  f(x) K(x,x')g(x')  -
 \int_{\substack{||x'-x||=\eps\\
x,x'\in D}}\, f(x) K(x,x')g(x')
 =  \int_{\substack{D\times D \\ \|x-x' \| \geq \eps}} \,d''f(x){\scriptstyle{\wedge}} K(x,x')g(x')
\end{eqnarray*}

\noindent Puisque $ K(x,x')g(x') =  L(x,x') g(x')$, nous avons :
\begin{eqnarray*}
 J_\eps &:=&\int_{\substack{||x'-x||=\eps\\  x,x'\in D}}\,
f(x) K(x,x')  g(x')={\big(N \, \Vol B(0,1)\big)}^{-1} \eps^{-N}\,\int_D \,  \big( \int_{\partial B(x',\eps)}\,f(x) \tilde{B} (x,x')\big) g(x') \\
&=&{\big(N \,\Vol B(0,1)\big)}^{-1} \eps^{-N}\,\int_D \,  \big(\int_{\partial B(0,\eps)}\,f(x' + w) \tilde{A} (w)\big) g(x')
\qquad \mbox{ où } \tilde{A} {\mbox{ est défini en (\ref{def-Atilde})}}
\end{eqnarray*}
Raisonnant  alors comme dans la preuve de (\ref{Omega-fonda}),
nous obtenons $ \lim_{\eps\rightarrow 0} J_\eps = \int_D \,  f(x')
g(x')$\,. Une régularisation standard permet de conclure.\qed

\medbreak\noindent {\sc notations.} Lorsque $n=1$ et $m=0$ (resp.
$n=0$ et $m=1$) nous noterons
\begin{equation}\label{K0K1}
K_0 =\Psi^* \Om_0 \; \mbox{ (resp. }
K_1 =\Psi^* \Om_1\,),
\end{equation}
\noindent  où $\Om_0$ et $\Om_1$ sont définis dans la proposition
\ref{Omega0} et en (\ref{Omega1})

\section{Quelques propriétés des fonctions qS-différentiables.}

Une application $f$ Fréchet-différentiable sur $D$ domaine de
$\Lambda_{0}^{n}$ à valeurs dans $\Lambda$ est S-différentiable si
et seulement $d''f\equiv 0$ sur $D$. Nous avons déjà vu , au troisième
paragraphe, que la situation n'est pas tout à fait la même lorsque les
variables sont dans $\Lambda_{1}$.

\medbreak\noindent
 {\sc Définition} :
 $f:\,D\subset  \Lambda_{0}^{n}\times\Lambda_{1}^{m}\rightarrow\Lambda$ est quasiment S-différentiable ({\it qS-différentiable} en abrégé)
 si elle est Fréchet-différentiable sur $D$ et telle que $d''f\equiv 0$ sur $D$. \\

 \noindent Si $f$ est S-différentiable, alors elle est qS-différentiable.\\
 Si $m=0$, $f$ est S-différentiable si et seulement si elle est qS-différentiable.\\

 On suppose toujours vérifiées les conditions ($A_0$) et ($A_1$) sur $\Lambda$.

 \subsection{ Une  formule de représentation intégrale sur
la frontière distinguée $\partial_0 P$ d'un polydisque.}

\begin{proposition}\label{weierstrass} Soit $D$ un ouvert de $\Lambda_0$ ou $\Lambda_1$ ; si une
suite ($f_k$) de fonctions  qS-différentiables sur $D$ converge
vers $f$ uniformément sur tout compact de $D$, alors $f$ est
qS-différentiable sur $D$.
\end{proposition}

\noindent Ce résultat se démontre de façon classique dès lors que l'on peut valider
une différentiation sous le signe somme. Et ce fait découle des estimations
$|D^\alpha K_0(x,y)|\lesssim ||x-y||^{-|\alpha
|-p}$ uniformément en $x \in G\subset\subset\Omega$ ouvert relativement compact
dans $D$ et $y\in D \setminus \overline{\Omega}.$ \\
L'estimation  analogue $|D^\alpha K_1(x,y)|\lesssim
||x-y||^{-|\alpha |-q+1}$ donne le résultat lorsque $D \subset
\Lambda_1$.

\medbreak\noindent

\begin{proposition}
Soit $P=\prod_{j=1}^{n+m} \Delta_j (a_j ; r_j )$ où chaque
$\Delta_j, j=1,\dots ,n$ (resp. $j=n+1,\dots ,m$) est un
polydisque de $\Lambda_0$ (resp. $\Lambda_1$) ouvert non vide. Si
$f$ est continue sur $\overline{P}$ et séparément
qS-différentiable sur $\Om \subset \Lambda_0 ^{\, n} \times
\Lambda_1 ^{\, m}$ avec $\overline{P} \subset \Om$, alors pour
tout $x=(y,\theta) =(y_{1} ,\dots ,y_{n},\theta_1 ,\dots ,\theta_m
)\in P$ :
\begin{equation} \label{cauchy-polyd}f(y,\theta) = \int _{\partial_0
P}\, f(w_{1} ,\dots ,w_{n} , \tau_1 ,\dots ,\tau_m ) \prod_{j=1}^n K_0 ( w_{j} ,y_{j} )\,.\prod_{j=n+1}^m K_1 (\tau_{j},\theta_{j}  )
\end{equation}
\end{proposition}

\noindent On déduit tout d'abord aisément de la proposition \ref{weierstrass} que
$y_j \mapsto f(a_1, ..., a_{j-1},y_j, \dots ,a_n ,b )$ est
qS-différentiable sur $\Delta_j$ dès que $a_j \in
\overline{\Delta_i}\,,\, \forall i\neq j$ et  $b_i \in
\overline{\Delta_i}\,, \forall i=1,..,m$ (et mutatis mutandis pour
$\theta_j \mapsto f(a , b_1 ,..., b_{j-1},\theta_j, \dots ,b_m)$). \\

 Nous pouvons écrire, en supposant $n\neq 0$ par exemple
$$   f(y_{1} ,\dots , y_{n},\theta_1,\dots ,\theta_m ) = \int_{\partial \Delta_1}\,
f(w_{1} ,y_2 ,\dots ,y_{n}, \theta_1 ,\dots ,\theta_m )  K_0 ( w_{1} , y_{1})\,.$$
  \noindent puisque  $y_{2}\mapsto f(w_{1},y_2 , \dots, y_{n}  , \theta )$
  reste qS-différentiable pour $w_{1} \in \partial \Delta_1$,  on
  peut itérer le processus, et appliquer le théorème de Fubini ($f$
  étant continue sur $\overline{P}$).\\

 \subsection{Harmonicité  ; un théorème de Hartogs de qS-différentiablité séparée.}

 Etablissons tout d'abord quelques propriétés
(harmonicité, principe du maximum) qui sont les analogues en
super-analyse de propriétés des fonctions holomorphes.

 \begin{proposition}\label{harmonicite}  Soit $D$ un domaine de  $\Lambda_{0}^{n}\times\Lambda_{1}^{m}$.
 Si $f : D \rightarrow\Lambda$  est  qS-différentiable sur $D$, alors elle est  harmonique
sur $D$, donc $\mathbb{R}$-analytique sur $D$.
 \end{proposition}

 \noindent Soit un polydisque $P=\prod_{i=1}^{n}\,P_{i}(a_{i},r_{i}).\prod_{j=1}^{m}\,P_{j}(b_{j},\varrho_{j})\subset\subset
 D$ de frontière distinguée
 $\partial_{0}P$.

 \noindent D'après (\ref{cauchy-polyd}), pour $x=(y,\theta)=(y_{1},...,y_{n},\theta_{1},...,\theta_{m})\in \frac{1}{2}P$\\
 $f(x)=f(y,\theta)=\int_{\partial_{0}P}\,f(z,\zeta)\prod_{j=1}^{n}K_{0}(z_{j},y_{j})\prod_{j=1}^{m}K_{1}(\zeta_{j},\theta_{j})$\\
 Pour tout $z_{j}\in \partial P_{j}(a_{j},r_{j})$, (resp. tout $\zeta_{j}\in\partial P_{j}(b_{j},\rho_{j})$), le noyau
 $K_{0}(z_{j},y_{j})$ (resp. $K_{1}(\zeta_{j},\theta_{j})$) est Fréchet-$C^{\infty}$ sur
 $\{y_{j}/\|y_{j}-a_{j}\|< r_{j}/{2}\}$ (resp. $\{\theta_{j}/\|\theta_{j}-b_{j}\|
 < \rho_{j}/{2}\}$) et l'on a, pour tous les multi-indices $\alpha_{j}$ et $\beta_{j}$,
 avec ces conditions sur $z_{j}$, $\zeta_{j}$, $y_{j}$ et $\theta_{j}$:\\
 $\|D_{y_{j}}^{\alpha_{j}}K_{0}(z_{j},y_{j})\|\,+\,\|D_{\theta_{j}}^{\beta_{j}}K_{1}(\zeta_{j},\theta_{j})\|\leq Cste(r_{j},\rho_{j}).$\\
On déduit alors de (\ref{cauchy-polyd}) que  $f$ est de classe
$C^{\infty}$ au sens de Fréchet sur $\frac{1}{2}P$. \\
 Nous avons
\begin{eqnarray*}
&& \forall j=1,...,n;\;\forall k=0,...,p;\;
\frac{\partial^{2}f}{\partial (y_{j}^{k})^{2}}\,=\,\frac{\partial}{\partial y_{j}^{k}}
\left(e_{k}\frac{\partial f}{\partial y_{j}^{0}}\right)\,=\,e_{k}^{\;2}\frac{\partial^{2}f}{\partial (y_{j}^{0})^{2}}\\
&& \forall j=1,...,m;\;\forall {\ell} = 1,...,r;
\forall k, s_{\ell} \leq k < s_{\ell +1} \;
\frac{\partial^{2}f}{\partial (\theta_{j}^{k})^{2}}\,=\,a_{k}^{\;2}\frac{\partial^{2}f}{\partial (\theta_{j}^{0})^{2}}
\end{eqnarray*}

Comme $\sum_{k=0}^{p}\,e_{k}^{\, 2}\,=\,\sum_{k=s_{1}}^{s_{2}
-1}\,a_{k}^{\,2}\, =\dots  = \sum_{k=s_{r} }^{q} a_k ^{\, 2}=0$, $f$
est harmonique.   \qed

\smallbreak\noindent Remarque : en utilisant l'hypoellipticité du
laplacien, nous aurions pu prouver l'harmonicité de $f$ sans avoir
recours à une représentation intégrale.

\vspace{2mm}

Une fonction $f$ séparément holomorphe par rapport à chaque
variable $z_j \,$ de $\C^n$ ,\,$ j=1,\dots , n$,  sur un domaine
$\Om$ de $\C^n$ est holomorphe sur $\Om$ sans autre hypothèse de
régularité globale sur $f$ mais ce théorème de Hartogs n'est plus
valable pour les fonctions  $\R$-analytiques. Nous déduisons
toutefois de la proposition \ref{harmonicite} et d'un résultat de
P. Lelong \cite{L} :

\begin{theoreme}\label{hartogs-separe}  Une
fonction $f$ définie sur un domaine $D$ de $\Lambda_0^{\,n} \times \Lambda_1 ^{\,m}$ à
valeurs dans $\Lambda$, séparément qS-différentiable
sur $D$ par rapport à
chacune de ses variables appartenant à $\Lambda_0$ ou $\Lambda_1$,
 est qS-différentiable sur $D$.
\end{theoreme}

\noindent La fonction $f$ est séparément harmonique en les
variables $y_j \in \Lambda_0 \approx \R^{p+1}$, $j=1,..,n$ et
$\theta_j \in \Lambda_1 \approx \R^{q}$, $j=1,..,m$ au voisinage
de tout point de $D$ donc est globalement harmonique en $(y_1 ,
\dots , y_n , \theta_1 ,\dots , \theta_m )$ d'après \cite{L} p.
561 ; par suite elle  est de classe $C^1$ au sens de Fréchet sur
$D$ d'où le résultat d'après la définition (\ref{d''global}) de $d''$. \qed

\begin{corollary} Soit $f$ une fonction qS-différentiable sur un domaine $D$ de $\Lambda_0 ^{\,n} \times \Lambda_1 ^{\,m}$ et à valeurs dans $\Lambda$. Si $f$ admet (en norme) un maximum local en un point de $D$, $f$ est constante sur $D$.
\end{corollary}

\noindent On déduit de l'harmonicité des composantes $f_j$ de
$f=\sum_0^{p+q}\, f_j e_j$ que la fonction continue $||f||^2$ est
sous-harmonique, donc localement constante dès qu'elle admet un
maximum local. D'où le résultat puisqu'alors $0 = \Delta ||f||^2
=2 \sum_0^{p+q} |\nabla f_j |^2$. \qed

\smallbreak\noindent {\sc  Remarque } : la présence de diviseurs de zéro interdit la généralisation de théorèmes de type
zéros isolés ou application ouverte.\\
Si l'on considère la CSA $\Lambda_0$ de l'exemple \ref{exemple-table} et la fonction $f$ S-différentiable  d'une variable $y$
appartenant à $\Lambda_0$ définie par $f(y) =y e_2$, l'ensemble des zéros de $f$ contient
entre autres $\{e_3 + \frac{e_2}{n}\,,\, n\in  \N^*\}$ et par ailleurs $f(\Lambda_0) =\mbox{ Vect } (e_2 ,e_3 )$ est non
ouvert dans $\Lambda_0$.\\

\subsection{Analyticité.}
Pour $j=1,...,m$ et $\theta_{j}=
 \sum_1^{q} \theta_j^i \eps_i $\, notons :
\begin{equation}\label{Z(theta)}
Z(\theta_j )=\,\sum_{i=1}^{q}\,a_{i}\theta_{j}^{i} \; \mbox{ et,
pour tout } k=1,\dots ,r \;:
  \,\, Z_k (\theta_j )=
\sum_{i=s_k}^{s_{k+1}
-1}\theta_j^i\, a_i
\end{equation}

\noindent Une fonction $f$ qS-différentiable sur $D \subset
\Lambda_{0}^{n}\times\Lambda_{1}^{m}$ est donc analytique réelle
en les variables $(y^0_{1},..., y^p_{1},\dots, $  $y^0_{n},\dots
,y^p_{n}, \theta^1_{1},...,\dots,\theta^q_1 ,\dots , \theta^1_{m}
,\dots , \theta^q_m)$. \\
En fait, une fonction qS-différentiable possède une propriété
d'analyticité bien plus forte puisqu'elle est localement
développable en série entière des variables $y_{1}=
\sum_{k=0}^p\,y^k_{1} e_k,\,y_2 , ...,y_{n},Z_1(\theta _1),...,
Z_{r}(\theta _ 1),\dots,$  $\dots , Z_1(\theta
_m),...,Z_{r}(\theta _m) $.

\begin{proposition}\label{S-analyticite}
Soit
$f:\,D\subset\Lambda_{0}^{n}\times\Lambda_{1}^{m}\rightarrow\Lambda$
une fonction qS-différentiable sur $D$. \\
Pour tout $(b,\beta)=(b_{1},...,b_{n},\beta_{1},...,\beta_{m})\in
D$, il existe $r>0$ et des scalaires $A_{I,J}\in\Lambda$ où I et
$J=(J_1 ,\dots ,J_r )$ sont des multi-indices de $\mathbb{N}^{n}$
et $(\mathbb{N}^{m})^r$ respectivement, tels que pour
$\|y_{j}-b_{j}\|<r,\;j=1,...,n \mbox{ et
}\;\;\|\theta_{j}-\beta_{j}\|<r,\;j=1,...,m$, nous ayons avec les
notations (\ref{les Z}) :
\begin{equation}\label{dev-analyticite}
f(y,\theta)\,=\,\sum_{I,J_{1} ,\dots , J_{r}}\,A_{I,J}(y-b)^{I}
 (Z_1  (\theta -\beta ))^{J_{1}}...(Z_r (\theta -\beta ))^{J_{r}}
 \end{equation}
avec absolue convergence de la série. De plus :
\begin{equation}\label{coeff-serie}
 I! J_1 ! \dots J_r ! \, A_{I,J} =
\dfrac{\partial ^{|I|+|J_1 |+\dots +|J_r |} \, f}{\partial y ^I \partial (Z_1 \theta)^{J_1}\dots \partial (Z_r \theta)^{J_r}}(b,\beta )
\end{equation}
Le développement (\ref{dev-analyticite}) est valable sur tout
polydisque $P$ de centre ($b,\beta $) relativement compact dans
$D$. Réciproquement, toute fonction de la forme
(\ref{dev-analyticite}) est qS-différentiable sur le domaine de
convergence de la série.
\end{proposition}

\medbreak\noindent Nous avons noté
\begin{equation}\label{les Z}
\begin{array}{l}
(y-b)^{I}=(y_{1}-b_{1})^{i_{1}}...(y_{n}-b_{n})^{i_{n}}\;si\,I=(i_{1},...,i_{n})\in\mathbb{N}^{n}\, ;\\
Z_k  (\theta -\beta ) = \Big( Z_k  ( \theta_{1}-\beta_{1}), \dots , Z_k ( \theta_{m}-\beta_{m}) \Big)   \mbox{ où }
Z_k (\theta_i ) \mbox{  est défini
en  }(\ref{Z(theta)}). \\ \mbox{ Si } J_k =(j_{k} ^1, \dots ,
j_{k}^{m} )\in \N^m \,,\\
 \Big( Z_k  (\theta -\beta ) \Big)^{J_{k}} =\Big( Z_k
(\theta_{1}-\beta_{1})\Big)^{\, j_k^1} \dots \Big( Z_k (
\theta_{m}-\beta_{m})\Big) ^{\, j_k ^{m}}  , \; \mbox{ et }
\dfrac{\partial ^{|J_k |}}{\partial (Z_k \theta)^{J_k}}
= \dfrac{\partial ^{|J_k |}}{\partial (\theta_1^{s_k })^{j_k^1}\dots \partial (\theta_m^{s_k })^{j_k^m} }    \;.
 \end{array}
 \end{equation}

\medbreak\noindent \begin{remark}\label{gadea} {\rm   Un théorème
de Gadea-Mu{\~n}oz (cf. \cite{GM}) établit  que pour une
$\R$-algèbre commutative $A$ de dimension finie, l'égalité entre
l'espace des fonctions $A$-différentiables sur $D$ ouvert non vide
de $A^n$ et celui des fonctions $A$-analytiques  sur $D$ (au sens
de notre proposition) est réalisée si et seulement si $A$ est une
$\C$-algèbre. Notre proposition \ref{S-analyticite} est-elle en
accord avec ce résultat de \cite{GM}  lorsque $m=0$ et donc,
$\Lambda = \Lambda_0$  algèbre commutative ? Oui, et nous prouvons
dans la section suivante que la condition $(A_0 )$ implique pour
$\Lambda_0$ d'être une $\C$-algèbre}.
\end{remark}

\begin{lemma}\label{poly-lambda0}
Soit $f : \Lambda_0 \rightarrow \Lambda$ un polynôme en les
variables réelles $(y^0 ,.., y^p )$ à coefficients dans $\Lambda$
;  $f$ est qS-différentiable si et seulement si $f$ est un
polynôme en la variable $y=\sum_{j=0}^p y^j e_j$.
\end{lemma}

\noindent Il suffit de  prouver la condition nécessaire pour les
polynômes homogènes car la partie homogène d'ordre k d'un polynôme
qS-différentiable est
qS-différentiable. On raisonne par récurrence sur le degré d.\\
Si $d=1$, $f$ est $\R$-linéaire (à une constante près) et S-différentiable donc de la
forme $f(y)=ay\;+\;f(0)$ où $a\in \Lambda$.\\
Soit $f$ polynôme homogène en les variables réelles $(y^0 ,.., y^p
)$ de degré $d\geq 2$ et S-différentiable. En tant que polynôme,
$f$ est lisse au sens de Fréchet, et l'on a pour tout $k=0,..,p$ :
$$ \forall j : \frac{\partial^2 f}{\partial y^j\partial y^k}=
\frac{\partial}{\partial y^k} (e_j \frac{\partial f}{\partial y^0}
)= e_j \frac{\partial}{\partial y^0} ( \frac{\partial f}{\partial
y^k} );$$ \noindent les $\frac{\partial f}{\partial y^k} $ sont
donc S-différentiables, homogènes de degré $d-1$, donc, par
récurrence, il existe $b\in \Lambda_0$ tel que $\frac{\partial
f}{\partial y^0} = b (y)^{d-1}$ et $\frac{\partial f}{\partial
y_k} = e_k b (y)^{d-1}$ ; par intégration $f(y)= \int_0^1 df(ty).y
\,\textrm{d}t= \frac{b}{d}(y)^d$.\qed

\begin{lemma}\label{poly-lambda1}
Soit $f : \Lambda_1 \rightarrow \Lambda$ un polynôme en les
variables réelles $(\theta^1 ,.., \theta^q )$ à coefficients dans
$\Lambda$ ; $f$ est qS-différentiable si et seulement si  $f$ est
un polynôme en les variables $Z_1(\theta ),...,Z_{r}(\theta )$.
\end{lemma}
\noindent La condition est clairement suffisante. On raisonne
encore par récurrence sur d pour la preuve de la nécessité de la
condition. Pour un polynôme homogène de degré $d=1$, le résultat
découle de la définition de la qS-différentiabilité . Il est aisé
de vérifier la qS-différentiablité de $\frac{\partial f}{\partial
\theta^{s_j}}$ pour $j=1,\dots ,r$ si $f$ est homogène de degré
$d\geq 2$ ; par hypothèse de récurrence
$$\frac{\partial f}{\partial \theta^{s_j}} (\theta)=
\sum_ {k_1 + \dots +k_r =d-1} \,b_K^{\,(j)}
\big(Z_1(\theta)\big)^{k_1}...\big(Z_{r}(\theta)\big)^{k_r}
\;,\; b_K^{\,j}\in \Lambda\,, j=1,\dots ,r \;;$$

\noindent D'où le résultat en appliquant la formule de Taylor et
les relations $ \frac{\partial f}{\partial \theta^{k}} = a_k
\frac{\partial f}{\partial \theta^{s_j}}   $ pour $ s_j \leq k <
s_{j+1}$, $k=1,\dots ,q$.   \qed

\begin{lemma}\label{poly-ttes-variables}
Soit $f : \Lambda_0 \times \Lambda_1 \rightarrow \Lambda$ un
polynôme de
degré d en les variables réelles $y^0 , ...,y^p , \theta^1 ,..., \theta^q$ ; \\
 $f$ est qS-différentiable si et seulement si   $f$ s'écrit sous la
 forme :
$$f(y,\theta)\,=\,\sum_{\underset{K=(k_0,k_1,\dots ,k_r )}{|K|=0}}^{d}\,A_{K}\, y^{k_0}
\big(Z_1 (\theta)\big)^{k_1}\dots \big(Z_r (\theta
)\big)^{k_r} \,\; \mbox{ où } A_{K}\in\Lambda,\;\forall K\in
\N^{r+1}\, .$$
\end{lemma}

\noindent Et nous avons un énoncé analogue pour un
polynôme  à valeurs dans $\Lambda$ en les variables réelles $y^0_1
, ...,y^p_1, \dots, $ $y^0_n , ...,y^p_n , \theta^1_1 ,...,
\theta^q_1 \dots \theta^1_m ,..., \theta^q_m$.

\noindent {La  preuve du lemme s'effectue par récurrence sur le
degré du polynôme et repose sur les deux lemmes précédents.}

\medbreak\noindent {\sc retour à  la proposition
\ref{S-analyticite}}: \\
Soulignons ici que, en dehors de la diagonale,
le noyau reproduisant $K_0 (x,y)$ des fonctions S-différentiables
sur un domaine de $\Lambda_0$ est lisse au sens de Fréchet, mais
n'est pas $S$-analytique ; nous ne pouvons donc pas raisonner
exactement comme dans le cas holomorphe de une ou plusieurs
variables, même pour passer du cas $n=1$ au cas $n$ quelconque. Et
la même remarque vaut pour le noyau $K_1$.

\smallbreak\noindent  -- La qS-différentiabilité (sur le
domaine de convergence de la série) d'une fonction admettant un
développement de type (\ref{dev-analyticite}) découle du lemme
\ref{poly-ttes-variables}.

\smallbreak\noindent  -- La preuve de (\ref{dev-analyticite}) et (\ref{coeff-serie}) sous l'hypothèse de qS-différentiabilité étant un peu longue et technique, nous la reportons en appendice.

\medbreak\noindent {\sc inégalités de type
  Cauchy }:\\
Lorsque $f$ est continue  sur $\bar P$ (ou seulement bornée en
norme sur $P$) et qS-différentiable sur $P= \prod_j \Delta_j (a_j
; r_j )$, elle est de classe $C^\infty $ au sens de Fréchet  sur
$P$ et on peut prendre les dérivées dans la formule
(\ref{cauchy-polyd}) pour $||y_j -a_j ||< r_j /2$ pour tout $j=1,
\dots , n$ et $||\theta_j -a_{n+j}||\leq r_{n+j} \,/2$ pour
$j=1,..,q$.

\noindent $\mbox{Rappelons que } \frac{\partial f}{\partial y_j} =
\frac{\partial f}{\partial y_j ^0 }  \mbox{ et }
  \frac{\partial ^{|\beta |}f}{\partial (Z_k (\theta ) ^\beta} =
  \frac{\partial ^{|\beta |}f}{\partial (\theta_1^{s_k}) ^{\beta_1} \dots \partial (\theta_m^{s_k}) ^{\beta_m}} \;.$
Nous obtenons :
\begin{proposition}\label{prop-inegalites-cauchy}
Pour $f$  qS-différentiable sur un polydisque $P=
\displaystyle{\prod_{j=1}^{n+m}} \Delta_j (a_j ; r_j )$ et
continue sur $\bar P$, nous avons :
\begin{equation}\label{ineg-cauchy}
 ||\frac{\partial^{|I |+|J_1|+\dots |J_r |} \, f}{\partial y^I \partial (Z_1 (\theta )) ^{J_1}\dots
\partial (Z_r (\theta )) ^{J_r}}   (a)|| \leq c_{I, J_1 ,...,J_r}\,  . \sup_{\bar P} ||f|| .  r^{-[I, J_1 ,...,J_r ]} \;
\end{equation}
$\mbox{ où } c_{I,J_1 ,...,J_r} \sim I ! J_1 !  \dots J_r !  $
lorsque $\min (|I|, |J_1 |,\dots ,|J_r | )\rightarrow +\infty$.
\end{proposition}

\noindent Nous avons noté  pour $I=(i_1 ,...,i_n )$ et $J_k =(j_1^k ,\dots , j_m^k )\;, k=1,...,r$ :\\
$r^{- [I,J_1 , \dots ,J_r ]}=r_1 ^{\,- i_1} \dots r_n ^{\,- i_n} (
r_{n+1}) ^{- j_1^1 -\dots - j_1 ^r}\dots (r_{n+m}) ^{ -j_m^1
-\dots - j_m ^r}$.

\smallbreak Nous  déduisons de ces inégalités et de la proposition
\ref{S-analyticite} un résultat de type Liouville :

\begin{corollary}\label{liouville}
{\it Une fonction qS-différentiable et bornée en norme sur
$\Lambda_0^{\,n}\times \Lambda_1^{\, m}$  est constante.}
\end{corollary}

\medbreak

\subsection{Un théorème de prolongement de type Hartogs-Bochner. }

Grâce à la formule de représentation intégrale des fonctions
qS-différentiables donnée par le corollaire \ref{corol-repres} et
aux propriétés de l'opérateur $d''$, nous obtenons :

\begin{theoreme}: Sous les conditions ($A_0$) et ($A_1$), si $\partial \Om$ est le bord connexe d'un domaine
$\Om$ borné de $\Lambda_0 ^{\,n} \times \Lambda_1 ^{\, m}$, avec
$n+m \geq 2$,  et $f$ une fonction qS-différentiable dans un
voisinage connexe de $\partial \Om$, alors $f$ se prolonge en une
fonction qS-différentiable sur $\Om$.
\end{theoreme}

L'opérateur $d''$ et le noyau $K$ vérifient hors de la diagonale de $\mathbb{R}_{\Lambda}^{n,m}$:
 $d''_x K^{(0)}(w,x)= - d''_w K^{(1)}(w,x)$ où $K^{(j)} (w,x)$ désigne la composante de
degré $j$ en $x$ et $n(p+1)+mq -1 -j$ en $w$ du noyau $K$. Par
suite, si $D_1$, $D_2$ sont des domaines à frontière de classe
$C^1$ par morceaux, $D_1 \subset \subset \Om \subset \subset D_2$
avec $\partial D_j \subset V$, $j=1,2$
 où $V$ est un voisinage ouvert de $\partial \Om$ sur lequel $f$
 est qS-différentiable, les fonctions $F_j = \int_{\partial D_j}\, f(w)
 K^{(0)} (w,x) $, $j=1,2$  sont qS-différentiables respectivement sur $\R^{n,m}_\Lambda \setminus \partial D_1$
 et sur $\R^{n,m}_\Lambda \setminus \partial D_2$. Supposant par exemple $n\neq
 0$, en appliquant le corollaire \ref{liouville}  à la fonction $(y_1
 ,a', b)\in \Lambda_0 \times \Lambda_0 ^{\, n-1}\times \Lambda_1
 ^m \mapsto F_1 (y_1 ,a' ;b)$ pour chaque $(a' ,b)$ tel que
 $\Lambda_0 \times \{(a' ,b)\} \cap \partial D_1 =\emptyset$, on
 montre ensuite, par prolongement analytique, que la fonction $F_1$ est nulle
 sur la composante connexe non bornée de $\overline{D_1}$. Par suite, la fonction
 $F_2 $ est qS-différentiable sur $\Om$ et prolonge $f$.
La preuve du théorème est alors analogue à celle du Kugelsatz de
Hartogs-Bochner donnée dans \cite{LM} par exemple.

\section{Commentaires sur les conditions $(A_ 0 )$ et $(A_1 )$.}

Nous avons déjà souligné, dès l'introduction, que la condition
$(A_0 )$ est nécessaire et suffisante si l'on veut obtenir une
solution fondamentale pour l'opérateur $d''$ opérant sur  les
formes différentielles  définies sur un ouvert de $\Lambda_0
^{\,n}$, $n\geq
1$. \\
La nécessité de conditions algébriques  pour obtenir une
représentation intégrale de
 fonctions satisfaisant une équation $\delta f=0$ pour un
 opérateur différentiel $\delta$ défini sur un certain type d'espaces peut paraître à première vue surprenante
 mais elle est "raisonnable" et dépend de l'opérateur étudié $\delta$.\\
 De telles conditions algébriques ont été, par exemple, mises en évidence par J. Ryan dans (\cite{Ry})
  ; dans cet article l'auteur obtient comme condition
 nécessaire et suffisante pour l'obtention d'une formule "généralisée" de représentation intégrale de type Cauchy  le fait que $A$ contienne une sous-algèbre complexe
 isomorphe à une algèbre de
 Clifford de dimension finie.\\
 Soulignons que tant dans \cite{GM} (cf. remarque \ref{gadea}) que dans \cite{Ry}, deux articles où des conditions de nature algébrique découlent de résultats d'analyse, les auteurs ne se placent pas dans le contexte de la super-anayse.

 \bigbreak   Que peut-on dire d'une super-algèbre $\Lambda$
 telle que $\Lambda_0$ vérifie $(A_0 )$ ? \\
 L'algèbre
 donnée au premier paragraphe par sa table (cf. exemple \ref{exemple-table}) apparaît comme une
 $\C$-algèbre (notons que dans ce cas $\Lambda_0 = \C \otimes_\R
 B_0$ où $B_0$ est la sous-algèbre réelle monogène (donc associative et commutative) engendrée par l'élément $b=e_4$ de polynôme minimal
 $X^3$, et $\Lambda =\C \otimes_\R B$, où $B$ est une sous-algèbre
 réelle de $\Lambda$
 admettant pour base comme $\R$-espace vectoriel $(e_0, e_2
 ,e_ 4 ,\eps_1 ,\eps_3 ,\eps_4)$).\\
En fait, on peut démontrer :

\begin{proposition}\label{C-superalgèbre}
 Toute super-algèbre $\Lambda$ réelle unitaire de dimension finie telle que la sous-algèbre $\Lambda_0$
vérifie la condition $(A_0 )$ possède une structure de
$\C$-algèbre compatible avec sa structure de $\R$-algèbre.
\end{proposition}

\noindent {\sc preuve.} Dans une CSA $\Lambda$, la sous-algèbre
$\Lambda_0$ est contenue dans le centre de $\Lambda$,
 par suite la proposition découle des  deux lemmes suivants :

\bigbreak \noindent
\begin{lemma}\label{A-C-algèbre} Toute algèbre  $A$ réelle commutative unitaire de dimension finie et vérifiant $(A_0 )$
possède un élément $a$ tel que $a^2 =-1$.
 \end{lemma}
 D'après le thm principal de structure de Wedderburn pour les algèbres associatives
 de dimensions finies (cf \cite{A} ou \cite{P}  par exemple), la $\R$-algèbre $A$ est une
 somme directe (de sev) $A= S \oplus \mathcal{N}$ où $\mathcal{N}$ est le radical de $A$ et $S$
 une sous-algèbre semi-simple de $A$ isomorphe à $A/\mathcal{N}$ ; $S$ étant commutative est isomorphe à un produit
 $\prod_{j=1}^s  K_j$ où $K_j =\R$ ou $\C$ ;
 le neutre multiplicatif $e$ de $A$  appartient à  $S$  ; en effet,  $e= \sigma + \nu$ avec $\sigma\in S$ et $\nu \in
 \mathcal{N}$ ; on déduit de l'idempotence de $e$,
l'égalité $\sigma^2 =\sigma$, puis
  $\nu ^2 =\nu$, d'où $\nu =0$ par nilpotence de $\nu$.  \\
  S'il existe une base $(e, e_1 , ..., e_p )$
  de $A$ telle que $e +\sum_{k=1}^p e_k ^{\,2}=0$ avec pour  tout  $k=1,...,p$, $e_k = b_k +\nu_k$ où
  $b_k \in S$ et  $\nu_k \in \mathcal{N}$, nous obtenons l'égalité $e+\sum_{k=1}^p b_k ^{\, 2} =0$,
  et cette dernière égalité est possible seulement si tous les $K_j$ sont  égaux à $\C$.\\
  Nous avons, en identifiant $S$ et $\prod_{j=1}^s  K_j$ pour simplifier les notations, $e=(1, \dots , 1) \in S$ ; il suffit de choisir
  $a =(i,\dots ,i)$ pour conclure.

\begin{lemma}\label{A-contenant-i} Si $A$  est une  algèbre réelle unitaire, d'unité $e$,
 de dimension finie possédant dans son centre un élément $\iota$ tel que $\iota ^2 = -e$, alors
  $A$ est une $\C$-algèbre.
 \end{lemma}

Le cas $dim_\R A=2$ étant immédiat,  supposons $A$ de dimension
$n>2$. On construit alors une suite strictement croissante de
sous-espaces $V_k$ stables par multiplication par $\ii$,
 avec $V_0 = \mbox{Vect }(e, \ii )$, ..., $V_k = \mbox{Vect }(e, \ii , b_1 , \ii b_1 ,  \dots ,b_k ,\ii b_k )$
 (en effet, si $z \notin V_k$, la famille
   $( e,\ii , b_1 ,\ii b_1, \dots, b_k , \ii b_k ,z, \ii z)$, est
   $\R$-libre).
Nécessairement $n$ est pair et si $n=2p+2$, $A=V_p$.\\
     Posant alors pour $ \lambda = \alpha +i \beta \in \C$ et $x\in A$ :
   $\lambda x = \alpha x + \beta \iota x$, $A$ peut être regardée comme
   une $\C$-algèbre, l'application $(x,y)\in A \mapsto xy$
   devenant $\C$-bilinéaire puisque $\iota$ commute avec tous les
   éléments de $A$. \qed

\bigbreak En ce qui concerne la condition $(A_1 )$, qui n'a pas
été utilisée dans la preuve de la proposition
\ref{C-superalgèbre}, signalons que le résultat d'analyticité
énoncé dans la proposition \ref{S-analyticite} est d'autant plus
intéressant que les $Z(\theta_j )$ regroupent un
plus grand nombre de variables réelles $\theta_j^k$. \\
Ainsi dans l'exemple \ref{exemple-table}, à condition de bien
choisir les "tranches" $\eps_{s_j} ,\dots , \eps_{s_{j+1}-1}$, $1
\leq j \leq r$, il n'apparaît pas  dans le développement
(\ref{dev-analyticite}) uniquement des $Z(\theta_j)$ faisant
intervenir seulement deux variables réelles. Plus précisément,
puisque $\eps_1 \eps_4 =e_2 \neq 0$, il n'existe
pas d'élément $\eta \in \Lambda_1$ tel que $\Lambda_1 = \eta \Lambda_0$ ; donc nous devons déterminer deux "tranches" minimum.\\
 La condition $(A_1 )$ est vérifiée puisque
$\Lambda_1 = \eps_1 {\mbox{Vect}} (e_0 , e_1 , e_4 , e_5 ) +
\eps_4  {\mbox{ Vect }}(e_0 , e_1  )$ ; posant $a_{s_1}=e_0$,
$a_{s_2}=e_1$, $a_{s_3}=e_4$ et $a_{s_4}=e_5$, nous avons alors
$Z(\theta_1 )= \sum_{k=1}^4\,
\theta_1^k a_{s_k}$. \\
Remarquons également que $\mbox{ Vect }(e_0 , e_1 , e_4 , e_5 )$
n'est pas une sous-algèbre de $\Lambda$ ; donc a fortiori les
"tranches" de $\Lambda_1$ ne sont pas nécessairement associées à
des $\C$-algèbres.

\section{Appendice.}
\subsection{Preuve du corollaire \ref{corol-repres} (b)}
  Toute forme différentielle à coefficients à valeurs dans
$\Lambda$ en les $N =n(p+1)+mq$  variables réelles $y_i^j$,
${1\leq i\leq n}$, $0\leq j\leq p$ et $\theta_i^\ell$, ${1\leq
i\leq m}$, ${1 \leq j\leq q}$ peut donc s'écrire de
 manière unique en fonction des 1-formes de l'ensemble\\
 $\mathfrak{F}  =\{dY_i \,, {{i=1,\dots ,n}}\;;\; dZ_k (\theta_i )\;,{{i=1,\dots ,m \;, k=1,\dots
 ,r}}\}$, que nous dirons 1-formes de bidegré (1,0) et des 1-formes
 dites de bidegré (0,1) appartenant à l'ensemble \\
 $\mathfrak{H}= \{dy^j_i \,,\, {\textstyle{i=1,\dots ,n\,,j=1,\dots ,p}}\, \;;\; d\theta^\ell_i \;,{{i=1,\dots
 ,m \,,\, \ell =1,\dots ,q\,,\, \ell \neq s_k\,,\,\forall
 k=1,..., r}}\,\}$.

 \noindent Une forme homogène $f$ de degré $d$ s'écrit $f=\sum_{|I|=0}^{n+mr}\,
 f^{I}$, où pour chaque multiindice $I$, la forme $f^{I}$ est de
 bidegré $(|I|,d-|I|)$ ; chaque composante $\sum_{|I|=k} f^{I}$ est une
 forme "pure" de bidegré $(k, d-k)$ ;\\
 Les formes $\lambda (\theta_i )= $ $dZ_1
 (\theta_i )d\theta_i ^2 \dots d\theta_i^{s_2 -1} dZ_2 (\theta_i
 )d\theta_i ^{s_2 +1 }
 \dots dZ_r (\theta_i )d\theta_i ^{s_r
 +1} \dots
d\theta_i ^q$ sont de bidegré $(r,q-r)$ et les formes $\omega (y_i
)= dY_i dy_i^1 \dots dy_i^p$
 sont de bidegré $(1,p)$.\\
 Notons que le noyau $K(x,x')$ est de bidegré total ({\small i.e.
 relativement à (x,x')}) $(n+mr , N-n-mr-1)$.

\noindent $K_{(\alpha,\beta )}$ désigne la composante de $K$ de
bidegré
   $(\alpha ,\beta)$ relativement à la variable $x'$ avec la
   convention  $K_{(\alpha,-1)} \equiv 0$.

\smallbreak Montrons en un premier temps la formule suivante :  si
$f$ est une forme de degré $d$ de classe $C^\infty$ dans $D$
continue ainsi que $df$ sur $\overline{D}$, nous avons pour tout
$x'\in D$ :
\begin{equation}\label{formule-repres-pour-d}
f(x')=\int_{\partial D}\,f(x) K(x,x')\,+\,(-1)^{d-1}\int_{D}\,df(x){\scriptstyle{\wedge}}K(x,x')\,
+ \,(-1)^{d} \,d \int_{D}\,f(x) K(x,x')
\end{equation}
\noindent Il suffit de tester $f$ contre une forme $g(x')$ de degré
$N-d$, de classe $C^\infty$ à support compact dans $D$. Calculons:
\begin{equation}\label{stokes}
\begin{split}
\; &\int_{\partial (D\times D)}g(x')f(x)K(x,x')\,\\
&=\,\int_{\partial D_{x'}\times D_{x}+(-1)^{N}D_{x'}\times \partial D_{x}}g(x')f(x)K(x,x')\,=\,(-1)^{N}\int_{D} g(x')\int_{\partial D}f(x)K(x,x')\\
&=\,\int_{D\times D}dg(x')f(x)K(x,x')+(-1)^{N-d}g(x')df(x)K(x,x')+(-1)^{N}g(x')f(x)[\Delta]\\
&=(-1)^{N-d+1}\int_{D} g(x')d_{x'}\int_{D}f(x)K(x,x')+(-1)^{N-d}\int_{D} g(x')\int_{D}df(x)K(x,x')
+(-1)^{N}\int_{D} g(x')f(x').
\end{split}
\end{equation}
En identifiant les membres de droite des seconde et dernière
lignes de (\ref{stokes}), on obtient
(\ref{formule-repres-pour-d}).

Si $f$ est une forme de bidegré $(a,b),$ en égalant dans
(\ref{formule-repres-pour-d}) les composantes de bidegré $(a,b)$,
nous obtenons la formule du (b) du corollaire.\qed

\subsection{preuve de la proposition \ref{S-analyticite}.}

\smallbreak\noindent   Si
$f:\,D\subset\Lambda_{0} \rightarrow\Lambda$ est qS-différentiable
sur $D$, elle est réelle analytique et pour tout $a\in D$, il
existe des réels $r_k >0$ tels que
\begin{equation} \label{R-analytique}
f(y)= \sum_{J=(j_0 ,.., j_p)\in \N^{p+1}} \,\alpha_J (y^0 -a^0 )^{j_0}\dots (y^p -a^p )^{j_p}
\mbox{ où } \alpha_J \in \Lambda
\end{equation}
avec absolue sommabilité pour $|y^k -a^k | <r_k$. Chaque
composante homogène de la série étant S-différentiable, d'après le
lemme \ref{poly-lambda0} il existe pour tout $k\in \N$ un élément
$c_k \in \Lambda$ tel que
\begin{equation*}\label{coeff}
  \sum_{\underset{J=(j_0 ,..., j_p )}{|J|=k}}\alpha_J (y^0 -a^0 )^{j_0}\dots (y^p -a^p )^{j_p} =c_k (y-a)^k \;;
  \end{equation*}
\noindent d'où par associativité des familles sommables :
\begin{equation}\label{S-analytique-1var}
f(y)= \sum_0^\infty c_k (y-a)^k
\end{equation}
\noindent avec  convergence en norme de la série pour $||y-a||< \min_j r_j$.\\
En dérivant terme à terme le second membre de
(\ref{R-analytique}), on obtient, puisque $e_0$ est le neutre de
$\Lambda$, en tenant compte de (\ref{coeff}) où nous avons $c_k = \alpha_{(k,0,...,0)}$ :
\begin{equation*}
\frac{\partial ^{k} f}{\partial ({y^0}) ^{\, k}}(a)= k! \, c_k\;\mbox{ soit encore } \frac{\partial ^{k} f}{\partial y ^{\; k}}(a)= k! \, c_k
\end{equation*}
\noindent Remarque : tout revient donc à dériver terme à terme (\ref{S-analytique-1var}).
\medbreak\noindent  $\bullet$ Soit $f$ qS-différentiable sur
$D\subset \Lambda_0^{\, n }$, $n\geq 2$ ;  pour alléger l'écriture
supposons $n=2$ et notons $w$ et $y$ les deux variables
appartenant à $\Lambda_0$. Au voisinage d'un point $(a,b)\in D$ :
\begin{eqnarray*}
f(w,y)&=& \sum_{\underset{I=(i_0 ,..,i_p ),J=(j_0 ,..,j_p )}{I, J\in \N^{p+1}}}
\,\gamma_{I,J} (w^0 -a^0)^{i_0}\dots
( w^p -a^p)^{i_p}  \, (y^0 -b^0)^{j_0}\dots
( y^p -b^p)^{j_p} \\
 &= &\sum_{J}
\Big[\sum_{I} \gamma_{I,J} (w^0 -a^0)^{i_0}\dots
( w^p -a^p)^{i_p} \Big] \,(y^0 -b^0)^{j_0}\dots
( y^p -b^p)^{j_p}
\end{eqnarray*}
 avec sommabilité pour
$||w-a||<r, ||y-b||<\rho$. \\
Par propriétés des séries entières, pour tout $J$ fixé, la famille
$\big ( \gamma_{I,J} (w^0 -a^0)^{i_0}\dots ( w^p -a^p)^{i_p}\big
)_I $ est sommable, de somme $\sum_{k=0}^\infty \big( \sum_{|I|=k}
\gamma_{I,J} (w^0 -a^0)^{i_0}\dots ( w^p -a^p)^{i_p} \big)$ ; $f$
est S-différentiable sur $D$, et donc la somme $\sum_I
\gamma_{I,J} (w^0 -a^0)^{i_0}\dots ( w^p -a^p)^{i_p}  $ est
S-différentiable sur $\{w\,,\, ||w-a|| <r\}$ ; appliquant le lemme
\ref{poly-lambda0} aux composantes homogènes de cette somme, nous
obtenons
\begin{equation*}
f(w,y) =\sum_{J\in \N^{p+1}}\,\Big[ \sum_{k=0}^\infty c_{k,J} (w-a)^k \Big](y^0 -b^0)^{j_0}\dots
( y^p -b^p)^{j_p}=
\sum_{J,k}\,c_{k,J} (w-a)^k \, (y^0 -b^0)^{j_0}\dots
( y^p -b^p)^{j_p}
\end{equation*}
$\mbox{avec } c_{k,J}= \gamma_{(k,0,...,0),J}\nonumber$.
\smallbreak\noindent Les familles $ \big( c_{k,J}  (y^0
-b^0)^{j_0}\dots ( y^p -b^p)^{j_p} (w-a)^k \big)_J$ sont
sommables, mais $(w-a)^k$  peut être diviseur de zéro lorsque
$k\neq 0$, et nous devons justifier, pour tout $k\in \N^*$, la
sommabilité de $ \big( c_{k,J}  (y^0 -b^0)^{j_0}\dots
( y^p -b^p)^{j_p} \big) _{J}$.\\
nous avons $(w-a)^k =(w^0 -a^0 )^k e_0 + \sum'$ ; la famille $
\Big (c_{k,J}  (y^0 -b^0)^{j_0}\dots ( y^p -b^p)^{j_p} (w^0
-a^0)^k e_0 \Big )_J$ est sommable et il en sera de même pour
$\Big( c_{k,J} (y^0 -b^0)^{j_0}\dots ( y^p -b^p)^{j_p}\Big)_{J} $
dès que $w^0 \neq a^0$ ce qui est toujours vérifié pour un $w$
proche de $a$ convenablement choisi. Les sommes $\sum_{J}  c_{k,J}
(y^0 -b^0)^{j_0}\dots ( y^p -b^p)^{j_p} $ étant S-différentiables,
on conclut grâce au lemme \ref{poly-lambda0} :
$$f(w,y)= \sum_{\ell,k=0}^\infty\,  A_{k,\ell} (w-a)^k (y-b)^\ell\; \;
\mbox{ avec } A_{k,\ell}= \gamma_{(k,0,...,0),(\ell,0)}$$
Et l'on obtient par dérivation :
\begin{equation*}\label{coeff-S-serie}
\ell ! \,k! \, A_{k,\ell}= \frac{\partial ^{\ell +k} f}{\partial (w^0 ) ^{\, k}\,\partial (y^0 )^{\, \ell}} (a,b)\;
\mbox{ soit encore } \ell ! \, k!\, A_{k,\ell}= \frac{\partial ^{\ell +k} f}{\partial w^{\, k} \partial y^{\, \ell}}(a,b)
\end{equation*}
\medbreak\noindent $\bullet$ Soit $f$ qS- différentiable sur un domaine $D \subset \Lambda_0 ^{\, n}\times \Lambda_1^{\, m}$ ; supposons $n=m=1$, le cas général étant analogue. \\
La proposition \ref{harmonicite} permet d'écrire au voisinage de $(a,\beta )\in
D$ :
\begin{eqnarray}\label{S-analyticite-mix}
f(y,\theta)&=& \sum_{{\substack{I=(i_0 ,..,i_p )\in \N^{p+1}\\J=(j_1 ,..,j_q ) \in \N^{q}}}}
\,\gamma_{I,J} (y^0 -a^0)^{i_0}\dots
( y^p -a^p)^{i_p}  \, (\theta ^1 - \beta^1 )^{j_1}\dots
( \theta^q -\beta^q)^{j_q} \\
  &= &\sum_{I}
\big(\sum_{J} \gamma_{I,J} (\theta ^1 - \beta^1 )^{j_1}\dots
( \theta^q -\beta^q)^{j_q} \big)\,
(y^0 -a^0)^{i_0}\dots
( y^p -a^p)^{i_p}  \nonumber
\end{eqnarray}
Pour tout multi-indice $I\in \N^{p+1}$, la série $\sum_J \,
\gamma_{I,J}(\theta ^1 - \beta^1 )^{j_1}\dots ( \theta^q
-\beta^q)^{j_q}$ est absolument sommable, de somme
$\sum_{k=0}^\infty \, \sum_{|J|=k}\, \gamma_{I,J} (\theta ^1 -
\beta^1 )^{j_1}\dots ( \theta^q -\beta^q)^{j_q}$ ; par
qS-différentiablité de chaque partie homogène dans cette série, on
déduit du lemme \ref{poly-lambda1} l'écriture suivante de $f$ au
voisinage de $(a,\beta )$ :
$$f(y,\theta)= \sum_I \Big(\sum_{k_1 ,..., k_r \in \N}\, \alpha _{I,(k_1 ,..,k_r )}
\big( Z_1 (\theta -\beta )\big)^{k_1} \dots \big(Z_r (\theta -\beta ) \big)^{k_r}\Big) \, (y^0 -a^0)^{i_0}\dots
( y^p -a^p)^{i_p}  \,$$
$\mbox{ où }
 \alpha_{I,(k_1 ,..,k_r )}= \gamma_{I, (j_{s_1},0,..,0,j_{s_2},0,..,0,j_{s_r},0,..,0)}$.\\
\noindent On vérifie, comme précédemment, la sommabilité de chaque
famille $(\alpha_{I,k_1 ,..,k_r} \,
 (y^0 -a^0)^{i_0}\dots
( y^p -a^p)^{i_p} )_I$ ; si $k_1 +\dots k_r \neq 0$ ; choisissons
$\theta = \beta + \sum_{\ell =1} ^{r}   \eta_\ell \eps_{s_\ell} $,
$\eta_\ell$ réel, $0<\eta_\ell<<1$ ; alors $Z_1 (\theta -\beta
))\big)^{k_1} \dots \big(Z_r (\theta -\beta ) )
 = \eta_1^{k_1} \dots \eta_r^{k_r} e_0$ est inversible  ; la famille $\big( \alpha _{I,k_1 ,..,k_r}
  \big( Z_1 (\theta -\beta )\big)^{k_1} \dots \big(Z_r (\theta -\beta ))\big)^{k_r}\Big) \, (y^0 -a^0)^{i_0}\dots
( y^p -a^p)^{i_p} \big)_{I=(i_0 ,\dots ,i_p )}$ étant sommable
pour $\eta_1 ,..., \eta_r$ suffisamment petits,  il en est de même
de $(\alpha_{I,(k_1 ,..,k_r )} \,
 (y^0 -a^0)^{i_0}\dots
( y^p -a^p)^{i_p} )_I$.  Chaque somme $\sum_I(\alpha_{I,(k_1
,..,k_r )} \,
 (y^0 -a^0)^{i_0}\dots
( y^p -a^p)^{i_p} )_I$ est S-différentiable, ainsi que ses
composantes polynomiales homogènes, il découle du lemme
\ref{poly-lambda0} le développement cherché :
 $$f(y,\theta )= \sum_{\mu=1}^\infty \sum_{\underset{K=(k_1 ,..,k_r )}{|K|=0}}^\infty \,
  A_{\mu,K} \, (y-a)^\mu \, Z_1 (\theta -\beta )\big)^{k_1} \dots \big(Z_r (\theta -\beta )\big)^{k_r}\;
 \mbox{ avec } A_{\mu,K}= \alpha_{(\mu ,0,..,0),K} \;.$$
\medbreak Dérivant terme à terme dans (\ref{S-analyticite-mix}),
on obtient ensuite par unicité du développement en série entière :
$$ \mu ! k_1 !...k_r !\, A_{\mu,(k_1 ,...,k_r )} =
\frac{\partial ^{\mu +k_1 + ...+k_r } f} {\partial (y^0 )^\mu
\partial (\theta^{s_1})^{k_1}\dots \partial
(\theta^{s_r})^{k_r}}(a,\beta ) = \frac{\partial ^{\mu +k_1 +
...+k_r } f}{\partial y ^\mu \partial (\theta^{s_1})^{k_1}\dots
\partial (\theta^{s_r})^{k_r}}(a,\beta )$$
\noindent $\bullet$ L'avant-dernière assertion de la proposition
découle de (\ref{coeff-serie}), de la proposition
\ref{prop-inegalites-cauchy} et du théorème de
prolongement des fonctions analytiques de variables réelles. \qed

\end{document}